\documentclass{tac}

\usepackage[all]{xy}
\usepackage{amsmath,amssymb}

\newcommand{\Ring}{\mathbf{Ring}}
\newcommand{\dLat}{\mathbf{dLat}}
\newcommand{\CV}{\mathbb{C}_{\Vmod}(\mathcal{E})}
\newcommand{\Shv}{\mathbf{Shv}}
\newcommand{\Rep}{\mathfrak{R}}
\newcommand{\V}{\textbf{V}}
\newcommand{\Vmod}{\mathbb{V}}

\newcommand{\Rex}{kernel determinig property}
\newcommand{\limit}{\underrightarrow{lim}}
\newcommand{\LH}[1]{\textbf{LH}/{#1}}

\input diagxy

\usepackage[colorlinks=true]{hyperref}
\hypersetup{allcolors=[rgb]{0.1,0.1,0.4}}

\author{William Zuluaga}

\address{Departamento de Matem\'aticas, Universidad Nacional de C\'ordoba \& CONICET
}

\title {A Pierce representation theorem for varieties with BFC}

\copyrightyear{2018}

\keywords{BFC, Central Elements, Sheaves}
\amsclass{00A00}

\eaddress{wijazubo@hotmail.com}

\newtheorem{thm}{Theorem}
\newtheorem{theorem}{Theorem} 
\newtheorem{lem}{Lemma}

\newtheoremrm{rem}{Remark}

\let\pf\proof
\let\epf\endproof
\mathrmdef{Hom}
\mathbfdef{Set}

\begin{document}

\maketitle
\begin{abstract}
 We generalize the Pierce representation theorem for (commutative) rings with unit to other algebraic categories with Definable Factor Congruences by using tools from topos theory. Of independent interest, we prove that an algebraic category with right existential definable factor congruences is coextensive if and only if has center stable by complements. 

\end{abstract}

\tableofcontents


\section{Introduction}\label{Introduction}

By a variety with $\vec{0}$ and $\vec{1}$ we understand a variety $\textbf{V}$ for which there are $0$-ary terms $0_{1}$ , ..., $0_{n}$ , $1_{1}$ , ..., $1_{n}$ such that $\textbf{V} \models \vec{0}\approx \vec{1}\rightarrow x\approx y$, where $\vec{0}=(0_{1}, ..., 0_{n})$ and $\vec{1}=(1_{1}, ..., 1_{n})$. If $\vec{a} \in A^{n}$ and $\vec{b} \in B^{n}$, we write $[\vec{a}, \vec{b}]$ for the n-uple $((a_{1} , b_{1} ), ..., (a_{n} , b_{n})) \in (A \times B)^{n}$. If $A\in \textbf{V}$ then we say that $\vec{e}=(e_{1}, ..., e_{n})\in A^{n}$ is a \emph{central element} of $A$ if there exists an isomorphism	$\tau: A\rightarrow A_{1}\times A_{2}$, such that $\tau(\vec{e})=[\vec{0}, \vec{1}]$. Also, we say that $\vec{e}$ and $\vec{f}$ are a \emph{pair of complementary central elements} of $A$ if there exists an isomorphism $\tau: A\rightarrow A_{1}\times A_{2}$ such that $\tau(\vec{e})=[\vec{0}, \vec{1}]$ and $\tau(\vec{f})=[\vec{1}, \vec{0}]$. As it is well known, the direct
product representations $ A\rightarrow A_{1}\times A_{2}$ of an algebra $A$ are closely related to the concept of factor congruence. A pair of congruences $(\theta, \delta)$ of an algebra $A$ is a pair of complementary factor congruences of $A$ if $\theta \cap \delta = \Delta$ and $\theta \circ \delta = \nabla$. In such case $\theta$ and $\delta$ are called \emph{factor congruences}. In most cases, the direct decompositions of an algebra are not unique; moreover, in general the pair $(\vec{e}, \vec{f})$ of complementary central elements does not determine the pair of complementary factor congruences  $(ker(\pi_{1}\tau ), ker(\pi_{2}\tau))$ where the $\pi_{i}'s$ are the canonical projections and $\tau$ is the isomorphism between $A$ and $A_{1}\times A_{2}$. We call such property the \emph{determining property} (DP).
\begin{itemize}
\item[(DP)] For every pair $(\vec{e}, \vec{f})$ of complementary central elements, there is a unique pair $(\theta, \delta)$ of complementary factor congruences such that, for every $i = 1, ..., n$
\begin{center}
\begin{tabular}{ccc}
$(e_{i} , 0_{i} ) \in \theta$ and $(e_{i} , 1_{i} ) \in \delta$ & and & $(f_{i} , 0_{i} ) \in \delta$ and $(f_{i} , 1_{i} ) \in \theta$ 
\end{tabular}
\end{center}
\end{itemize}

Observe that (DP) is in some sense the most general condition guaranteeing that central elements have all the information about direct product decompositions in the variety. In \cite{SV2009} it was proved that (DP) is equivalent to each one of the following conditions:

\begin{itemize}
\item[(DFC)] $\textbf{V}$ has definable factor congruences; i.e, there is a first order formula $\psi(\vec{z} , x, y)$ such that for every $A, B \in \textbf{V}$
\begin{center}
$A\times B\models \psi([\vec{0}, \vec{1}], (a,b),(a',b'))$ iff $a=a'$
\end{center}
\item[(BFC)] $\textbf{V}$ has Boolean factor congruences, i.e., the set of factor congruences of any algebra in $\mathcal{V}$  is a Boolean sublattice of its congruence lattice.
\end{itemize}

Let $\textbf{V}$ a variety with BFC. If the formula $\psi$ of (DFC) is existencial we will say that $\textbf{V}$ is a variety with exDFC. The aim of this work is to exhibit a representation theorem for varieties with exDFC in terms of internal conected models in toposes of sheaves over a Bolean algebra. The present work is motivated by the Pierce's representation theorem for integral rigs \cite{Z2016} and Lawvere's strategic ideas about the topos-theoretic analysis of coextensive algebraic categories \cite{L2008}.

\section{Preliminaries}\label{Preliminaries}

\subsection{Notation and basic results}\label{Notation and basic results}

If $A$ is an algebra, we denote the congruence lattice of $A$ by $Con(A)$. As usual, the join operation of $Con(A)$ is denoted by $\vee$. If $f:A\rightarrow B$ is an homomorphism we write $Ker(f)$ for the congruence of $A$, defined by $\{(a,b)\in A\times A\mid f(a)=f(b)\}$. The universal congruence on $A$ is denoted by $\nabla^{A}$ and $\Delta^{A}$ denotes the identity congruence on A (or simply $\nabla$ and $\Delta$ when the context is clear). If $S\subseteq A$, we write $\theta^{A}(S)$ for the least congruence containing $S\times S$. If $\vec{a}, \vec{b} \in A^{n}$,
then $\theta^{A}(\vec{a}, \vec{b})$ denotes the congruence generated by $C=\{(a_{k}, b_{k}) \mid 1 \leq k \leq n\}$. If $\vec{a}, \vec{b} \in A^{n}$ and $\theta \in Con(A)$, we write $\vec{a} \equiv \vec{b}(\theta)$ or $[\vec{a},\vec{b}]\in \theta$ to express that $(a_{i} , b_{i}) \in \theta$, for
$i = 1,..., n$. We use $FC(A)$ to denote the set of factor congruences of A. A
variety $\textbf{V}$ has Boolean factor congruences if for every $A \in \textbf{V}$, the set $FC(A)$ is a distributive sublattice of $Con(A)$. We write $\theta \diamond \delta$ in $Con(A)$ to denote that $\theta$ and $\delta$ are complementary factor congruences of $A$. If $\theta \in FC(A)$, we use $\theta^{\star}$ to denote the factor complement of $\theta$. If $\theta, \delta\in Con(A)$ we say that $\theta$ and $\delta$ \emph{permutes} if $\theta \circ \delta= \delta \circ \theta$.
\\

A \emph{system over} $Con(A)$ is a $2n$-ple $(\theta_{1},...,\theta_{n}, x_{1},...,x_{n})$ such that $(x_{i},x_{j})\in \theta_{i}\vee\theta_{j}$, for every $i,j$. A \emph{solution} of the system $(\theta_{1},...,\theta_{n}, x_{1},...,x_{n})$ is an element $x\in A$ such that $(x,x_{i})\in \theta_{i}$ for every $1\leq i\leq n$. Observe that if $\theta_{1}\cap ... \cap \theta_{n}=\Delta^{A}$, thus the system $(\theta_{1},...,\theta_{n}, x_{1},...,x_{n})$ has at least one solution.

\lem\label{basics about systems} 
Let $\theta$ and $\delta$ be congruences of $A$. The following are equivalent:
\begin{enumerate}
\item $\theta$ and $\delta$ permutes.
\item $\theta \vee \delta =\theta \circ \delta$
\item For every $x,y\in A$, the system $(\theta, \delta, x,y)$ has a solution.
\end{enumerate}
\endlem

Given two sets $A_{1}$,$A_{2}$ and a relation $\delta$ in $A_{1} \times A_{2}$,  we say that $\delta$ \emph{factorizes} if there exist sets $\delta_{1} \subseteq A_{1} \times A_{1}$ and $\delta_{2} \subseteq A_{2} \times A_{2}$ such that $\delta=\delta_{1}\times \delta_{2}$, where \[\delta_{1}\times \delta_{2}=\{((a,b),(c,d))\mid (a,c)\in \delta_{1}, (b,d)\in \delta_{2}\}\]
So, if $\delta\in Con(A_{1}\times A_{2})$ factorizes in $\delta_{1}, \delta_{2}$ it follows that $\delta_{i}\in Con (A_{i})$, for $i=1,2$.
\lem[\cite{BB1990}]\label{FC Factors equivalent BFC}
Let $\V$ be a variety. The following are equivalent:
\begin{enumerate}
\item $\V$ has BFC.
\item $\V$ has factorable factor congruences. I.e. If $A,B\in \V$ and $\theta\in FC(A\times B)$, then $\theta$ factorizes.
\end{enumerate}
\endlem

We say that a variety has the \emph{Fraser-Horn property} (see \cite{FH1970}) (FHP) if every congruence on a (finite) direct product of algebras factorizes. 
\\

Given a variety $\V$ and a set of variables $X$, we use $\textbf{F}_{\V}(X)$ to denote the free algebra of $\V$ freely generated by $X$ (or simply $\textbf{F}(X)$, if the context is clear). If $X = \{x_{1} , . . . , x_{n} \}$, then we use $\textbf{F}_{\V}(x_{1} , . . . , x_{n})$ instead of $\textbf{F}_{\V} (\{x_{1} , . . . , x_{n} \})$. 
\\
As a final remark we should recall that all the algebras considered along this work will always have finite $n$-ary function symbols and its type (unless necessary) will be omitted.

\subsection{Generalities about Varieties with DFC}\label{Generalities about Varieties with DFC}


Let $\textbf{V}$ a variety with $\vec{0}$ and $\vec{1}$ and suppose that has DFC. For every $A\in \textbf{V}$, we write $Z(A)$ to denote the set of central elements of $A$ and $\vec{e}\diamond_{A} \vec{f}$ to denote that $\vec{e}$ and $\vec{f}$ are complementary central elements of $A$. If $\vec{e}$ is a central element of $A$ we write $\theta_{\vec{0}, \vec{e}}^{A}$ and $\theta_{\vec{1}, \vec{e}}^{A}$ for the unique pair of complementary factor congruences satisfying $\vec{e}\equiv \vec{0} (\theta_{\vec{0}, \vec{e}}^{A})$ and  $\vec{e}\equiv \vec{1} (\theta_{\vec{1}, \vec{e}}^{A})$. It follows that $\vec{0}$ and $\vec{1}$ are central elements in every algebra $A$ and the factor congruences associated to them are $\theta_{\vec{0}, \vec{0}}^{A}=\Delta^{A}$, $\theta_{\vec{1}, \vec{0}}^{A}=\nabla^{A}$ and $\theta_{\vec{0}, \vec{1}}^{A}=\nabla^{A}$, $\theta_{\vec{1}, \vec{1}}^{A}=\Delta^{A}$, respectively. If there is no place to confusion, we write $\theta_{\vec{0}, \vec{e}}^{A}$ and $\theta_{\vec{1}, \vec{e}}^{A}$ simply as $\theta_{\vec{0}, \vec{e}}$ and $\theta_{\vec{1}, \vec{e}}$. Since $\textbf{V}$ has BFC, factor complements are unique so we obtain the following fundamental result

\thm\label{Bijection betwen centrals and factor congruences}
Let $\textbf{V}$ a variety with DFC. The map $g:Z(A)\rightarrow FC(A),$ defined by $g(e)=\theta_{\vec{0}, \vec{e}}^{A}$ is a bijection and its inverse $h:FC(A)\rightarrow Z(A)$ is defined by $h(\theta)=\vec{e}$, where $\vec{e}$ is the only $\vec{e}\in A^{n}$ such that $\vec{e}\equiv \vec{0}(\theta)$ and $\vec{e}\equiv \vec{1}(\theta^{\ast})$.
\endthm

As a consequence of Lemma \ref{basics about systems}, the we obtain the following result for varieties with BFC 

\lem\label{BFC has permutable congruences} 
In every algebra $A$ of a variety $\textbf{V}$ with BFC every pair of factor congruences permutes. 
\endlem

Those facts, allows us to define some operations in $Z(A)$ as follows: Given $\vec{e}\in Z(A)$, the \emph{complement $\vec{e}^{c_{A}}$} of $\vec{e}$, is the only solution to the equations $\vec{z}\equiv \vec{1}(\theta_{\vec{0},\vec{e}})$ and $\vec{z}\equiv \vec{0}(\theta_{\vec{1},\vec{e}})$. Given $\vec{e}, \vec{f}\in Z(A)$, the \emph{infimum} $\vec{e}\wedge_{A}\vec{f}$ is the only solution to the equations $\vec{z}\equiv \vec{0}(\theta_{\vec{0},\vec{e}}\cap \theta_{\vec{0},\vec{f}})$ and $\vec{z}\equiv \vec{1}(\theta_{\vec{1},\vec{e}}\vee \theta_{\vec{1},\vec{f}})$. Finally, the \emph{supremum} $\vec{e}\vee_{A}\vec{f}$ is the only solution to the equations $\vec{z}\equiv \vec{0}(\theta_{\vec{0},\vec{e}}\vee \theta_{\vec{0},\vec{f}})$ and $\vec{z}\equiv \vec{1}(\theta_{\vec{1},\vec{e}}\cap \theta_{\vec{1},\vec{f}})$. 
\\

As result, we obtain that $\textbf{Z}(A)=(Z(A),\wedge_{A},\vee_{A}, ^{c_{A}},\vec{0},\vec{1})$ is a Boolean algebra which is isomorphic to $(FC(A), \vee, \cap, ^{\ast},\Delta^{A},\nabla^{A})$. Also notice that $\vec{e}\leq_{A} \vec{f}$ iff $\theta^{A}_{\vec{0},\vec{e}}\subseteq \theta^{A}_{\vec{0},\vec{f}}$ iff $\theta^{A}_{\vec{1},\vec{f}}\subseteq \theta^{A}_{\vec{1},\vec{e}}$. If the context is clear enough, we will not use the subscripts in the operations of $\textbf{Z}(A)$. 
\\

We conclude this section with a result which will be useful in Section \ref{the representation associated to a Vmodel}.

\lem[\cite{B2012}]\label{Useful lema Centrals}
Let $\V$ be a variety with DFC, $A\in \V$. For every $\vec{e},\vec{f}\in Z(A)$, the following holds:
\begin{enumerate}
\item $\vec{a}=\vec{e}\wedge_{A}\vec{f}$ if and only if $[\vec{0},\vec{a}]\in \theta_{\vec{0},\vec{e}}$ and $[\vec{a},\vec{f}]\in \theta_{\vec{1},\vec{e}}$.
\item $\vec{a}=\vec{e}\vee_{A}\vec{f}$ if and only if $[\vec{1},\vec{a}]\in \theta_{\vec{1},\vec{e}}$ and $[\vec{a},\vec{f}]\in \theta_{\vec{0},\vec{e}}$.
\end{enumerate}
\endlem

\section{The universal property}\label{The universal property}

In the Introduction we saw that for every variety with $\vec{0}$ and $\vec{1}$, having BFC is equivalent to the variety having definable factor congruences. In this section we introduce several definitions concerning with the different sorts of definability that arise at the light of this context. In addition, we present some useful results that arise from the universal property of principal congruences in varieties with BFC.

\definition\label{Definability definitions}
Let $\textbf{V}$ a variety with BFC. 
\begin{enumerate}
\item A formula $\rho(\vec{z},x,y)$ defines $\theta_{\vec{1}, \vec{e}}$ in terms of $\vec{e}$ if for every $A,B\in \textbf{V}$, $a,b\in A$ and $c,d\in B$ 
\[A\times B\models \rho([\vec{0},\vec{1}], (a,c), (b,d))\;\textrm{iff}\; c=d \] 
\item A formula $\lambda(\vec{z},x,y)$ defines $\theta_{\vec{0}, \vec{e}}$ in terms of $\vec{e}$ if for every $A,B\in \textbf{V}$, $a,b\in A$ and $c,d\in B$ 
\[A\times B\models \lambda([\vec{0},\vec{1}], (a,c), (b,d))\;\textrm{iff}\; a=b \]
\end{enumerate} 
\enddefinition

In the last case, we also say that $\rho$ \emph{defines $\theta_{\vec{0}, \vec{e}}$ in terms of $\vec{e}^{c}$}.
\\

Notice that if a formula $\rho$ defines $\theta_{\vec{1}, \vec{e}}$ in terms of $\vec{e}$, for every algebra $A\in \textbf{V}$ and $\vec{e}\in Z(A)$, it follows that $\theta^{A}_{\vec{1},\vec{e}}=\{(a,b)\mid A\models \rho(\vec{1},a,b)\}$. A similar statement is obtained when a formula $\lambda$ defines $\theta_{\vec{0}, \vec{e}}$ in terms of $\vec{e}$.
\\

Altough in \cite{SV2009}, it was proved that the items $1.$ and $2.$ of the Definition \ref{Definability definitions} are equivalent (which is not trivial, since in general $\vec{0}$ and $\vec{1}$ are not interchangeables), such equivalence does not preserve the complexity of the formulas (c.f. \cite{BV2013}). This situation motivates the need of introducing several definitions in terms of the complexity of the formulas envolved.    
\\

We say that a variety $\textbf{V}$ with $\vec{0}$ and $\vec{1}$ has \emph{right existentially defined factor congruences} (RexDFC) if the formula that defines $\theta_{\vec{1}, \vec{e}}$ in terms of $\vec{e}$ is existential. Analogously, if the formula that defines $\theta_{\vec{0}, \vec{e}}$ in terms of $\vec{e}$ is existential, we say that $\textbf{V}$ has \emph{left existentially defined factor congruences} (LexDFC). If $\textbf{V}$ has RexDFC and LexDFC, we say that $\textbf{V}$ has \emph{twice existentially defined factor congruences} (TexDFC). Similar definitions arise when the considered formula is positive or equational (a finite conjunction of equations). In the positive case, we use the acronyms RpDFC, LpDFC and TpDFC to mean that the variety has \emph{right positively defined factor congruences}, \emph{left positively defined factor congruences} and \emph{twice positively defined factor congruences}, respectively. For a further reading about  varieties with equationally definable factor congruences the reader can consult \cite{BV2013} and \cite{BV2017}.
\\

\lem[\cite{S2010}]\label{Existential implies Positive}
For every variety $\textbf{V}$ with BFC the following holds:
\begin{enumerate}
\item RexDFC implies RpDFC.
\item LexDFC implies LpDFC.
\item TexDFC implies TpDFC.
\end{enumerate} 
\endlem

The following result expose the intimate relation between $\theta_{\vec{0}, \vec{e}}$, and the complexity of the formula that defines it. 

\lem[\cite{B2012}]\label{Positive implies principal congruence}
Let $\textbf{V}$ be variety with BFC. 
\begin{enumerate}
\item If $\textbf{V}$ has RpDFC, then for every $A\in \textbf{V}$ and $\vec{e}$ central element of $A$ we get that $\theta_{\vec{1}, \vec{e}}^{A}=\theta^{A}(\vec{1},\vec{e})$.
\item If $\textbf{V}$ has LpDFC, then for every $A\in \textbf{V}$ and $\vec{e}$ central element of $A$ we get that $\theta_{\vec{0}, \vec{e}}^{A}=\theta^{A}(\vec{0},\vec{e})$.
\item If $\textbf{V}$ has TpDFC, then for every $A\in \textbf{V}$ and $\vec{e}$ central element of $A$ we get that $\theta_{\vec{1}, \vec{e}}^{A}=\theta^{A}(\vec{1},\vec{e})$ and $\theta_{\vec{0}, \vec{e}}^{A}=\theta^{A}(\vec{0},\vec{e})$.
\end{enumerate}
\endlem

We will say that an homomorphism $f:A \rightarrow P$ has the \emph{universal property of identify the elements of $S$}, if for every homomorphism $g:A\rightarrow C$, such that $g(a)=g(b)$, for every $a,b\in S$; there exists a unique homomorphism $h:B\rightarrow C$, such that the diagram

\begin{displaymath}
\xymatrix{
A \ar[r]^-{f} \ar[dr]_-{g} & B \ar@{-->}[d]^-{h}
\\
 & C
}
\end{displaymath}
\noindent
commutes.
\\

The following Lemma is an standard result in universal algebra. Nevertheless it is a key observation which will be useful for the rest of this paper. 

\lem\label{universal property principal congruences} Let $A$ be an algebra with finite $n$-ary function symbols and $S\subseteq A$. Then, the canonic homomorphism $\nu_{S}:A\rightarrow A/\theta(S)$ has the universal property of identify all the elements of $S$. 
\endlem

Recall that, as a consequence of Lemma \ref{universal property principal congruences}, we get that, for every $\vec{a},\vec{b}\in A^{n}$, the canonical homomorphism $A\rightarrow A/\theta(\vec{a},\vec{b})$ has the universal property of identify the elements of the set $[\vec{a},\vec{b}]$.

\corollary\label{corollary universal property} 
Let $\textbf{V}$ a variety with BFC. Then:
\begin{enumerate}
\item  If $\textbf{V}$ has RexDFC, for every $A\in \textbf{V}$ and $\vec{e}$ central element of $A$, the canonical morphism $A\rightarrow A/\theta(\vec{1},\vec{e})$ has the universal property of identify $\vec{e}$ with $\vec{1}$.
\item  If $\textbf{V}$ has LexDFC, for every $A\in \textbf{V}$ and $\vec{e}$ central element of $A$, the canonical morphism $A\rightarrow A/\theta(\vec{0},\vec{e})$ has the universal property of identify $\vec{e}$ with $\vec{0}$.
\end{enumerate}
\endcorollary
\pf
Apply Lemmas \ref{Existential implies Positive}, \ref{Positive implies principal congruence} and \ref{universal property principal congruences}. 
\epf

\lem\label{diagram pushout}
Let $A$ and $B$ be algebras with finite $n$-ary function symbols and $f:A\rightarrow B$ an homomorphism. Then, for every $S\subseteq A$, the diagram
\begin{displaymath}
\xymatrix{
A \ar[r]^-{\nu_{S}} \ar[d]_-{f} & A/\theta^{A}(S) \ar[d]
\\
B \ar[r]_-{\nu_{f(S)}} & B/\theta^{B}(f(S))
}
\end{displaymath}
\noindent
is a pushout. \endlem

\pf

Let $a,b\in S$ and consider the following diagram:
\begin{displaymath}
\xymatrix{
A \ar[r]^-{h_{S}} \ar[d]_-{f} & A/\theta^{A}(S) \ar[d]_-{k} \ar@/^1pc/[ddr]^-{\alpha} & 
\\
B \ar[r]^-{h_{f(S)}} \ar@/_1pc/[drr]_-{\beta} & B/\theta^{B}(f(S)) \ar@{-->}[dr]_-{\gamma} &
\\
 & & C
}
\end{displaymath}

Observe that,
\[h_{f(S)}(f(a))=h_{f(S)}(f(b)). \]

Then, by Lemma \ref{universal property principal congruences}, there exists a unique $k:A/\theta^{A}(S)\rightarrow B/\theta^{B}(f(S))$, such that the inner square commutes. Suppose now that $\alpha h_{S}=\beta f$. Thus, for $a,b\in S$ given, since $(a,b)\in \theta^{A}(S)$, we have that 
\[\beta(f(a))=\alpha (h_{S}(a)))=\alpha(h_{S}(b))=\beta(f(b)) \]

so again by Lemma \ref{universal property principal congruences}, there exists a unique $\gamma:B\rightarrow C$, suh that the downward triangle commutes. Finally, to verify that the upper triangle commutes, notice that 
\[ (\gamma k)h_{S}=\gamma (kh_{S})=\gamma (h_{f(S)}f)=\beta f= \alpha h_{S}\] 

Since $h_{S}$ is epi, we conclude that $\gamma k=\alpha$.
\epf

\corollary\label{corollary diagram pushout}
Let $\textbf{V}$ be a variety with BFC, $A,B\in \textbf{V}$, $f:A\rightarrow B$ be an homomorphism, $\vec{e}$ be a central element of $A$ and $f(\vec{e})=(f(e_{1}),...,f(e_{n}))$, then:
\begin{enumerate}
\item If $\textbf{V}$ has RexDFC, the diagram

\begin{displaymath}
\xymatrix{
A \ar[r]^-{\nu_{e}} \ar[d]_-{f} & A/\theta^{A}(\vec{1},\vec{e}) \ar[d]
\\
B \ar[r]_-{\nu_{f(e)}} & B/\theta^{B}(\vec{1},f(\vec{e}))
}
\end{displaymath}

is a pushout. 

\item If $\textbf{V}$ has LexDFC, the diagram

\begin{displaymath}
\xymatrix{
A \ar[r]^-{\nu_{e}} \ar[d]_-{f} & A/\theta^{A}(\vec{0},\vec{e}) \ar[d]
\\
B \ar[r]_-{\nu_{f(e)}} & B/\theta^{B}(\vec{0},f(\vec{e}))
}
\end{displaymath}

is a pushout. 
\end{enumerate}
\endcorollary
\pf
Apply Lemmas \ref{Existential implies Positive}, \ref{Positive implies principal congruence} and \ref{diagram pushout}.
\epf

\section{Coextensivity and Center Stability}\label{Coextensivity and Center Stabitlity}

In the context of varieties with BFC one may be tempted to think that in general, homomorphisms preserves central elements and even complementary central elements. Unfortunately that is not case. Even in varieties with BFC having good properties like the Frasier Horn ones, the preservation of central elements is restricted to surjective homomorphisms (c.f. \cite{V1996}). In this section we prove that the coextensivity of algebraic categories associated to varieties with RexDFC and center stable is equivalent to ask the variety having center stable by complements.

\definition\label{Coextensivity} A category with finite limits $\mathcal{C}$ is called extensive if has finite coproducts and the canonical functors $1 \rightarrow \mathcal{C}/0$ and $\mathcal{C}/X \times \mathcal{C}/Y \rightarrow \mathcal{C}/(X + Y )$ are equivalences. 
\enddefinition

If the opposite $\mathcal{C}^{op}$ of a category $\mathcal{C}$ is extensive, we will say that $\mathcal{C}$ is coextensive. Classical examples of coextensive categories are the categories ${\Ring}$ and $\dLat$ of commutative rings with unit and bounded distributive lattices. In the following, we will use a characterization proved in \cite{CW1993}.

\proposition\label{Coextensivity Proposition} 
A category $\mathcal{C}$ with finite coproducts and pullbacks along its injections is extensive if and only if the following holds:
\begin{enumerate}
\item (Coproducts are disjoint.) For every $X$ and $Y$, the square below is a pullback
\begin{displaymath}
\xymatrix{
0 \ar[r]^-{!} \ar[d]^-{!} & Y \ar[d]^-{in_{1}} 
\\
X \ar[r]_-{in_{0}}  & X+Y
}
\end{displaymath}

\item (Coproducts are universal.) For every $X,X_{i}, Y_{i}$ with $i=0,1$ and $f : X \rightarrow Y_{0}+Y_{1}$, if the squares below are pullbacks
\begin{displaymath}
\xymatrix{
X_{0} \ar[r]^-{x_{0}} \ar[d]_-{h_{0}} & X \ar[d]^-{f} & X_{1} \ar[l]_-{x_{1}} \ar[d]^-{h_{1}}
\\
Y_{0} \ar[r]_-{y_{0}} & Y_{0}+Y_{1} & Y_{1} \ar[l]^-{y_{1}}
}
\end{displaymath}
\noindent
then the cospan $X_{0}\rightarrow X \leftarrow X_{1}$ is a coproduct.
\end{enumerate}
\endproposition 

Let $\textbf{V}$ be a variety with BFC, $A,B\in \textbf{V}$ and $f:A\rightarrow B$ an homomorphism. We will say that $f$ \emph{preserves pairs of complementary central elements} if preserves central elements; i.e, for all $e\in Z(A)$ it follows that $f(e)\in Z(B)$ and furthermore, 
\[e_{1}\diamond_{A}e_{2} \Rightarrow f(e_{1})\diamond_{B}f(e_{2}) \]

If every homomorphism between the algebras of $\textbf{V}$ preserves central elements, we say that $\textbf{V}$ has \emph{stable center} (SC). If $\textbf{V}$ has SC and every homomorphism between the algebras of $\textbf{V}$ preserves central elements, we say that $\textbf{V}$ has \emph{center stable by complements} (CSC).

\rem\label{SC and CSC are not trivial}
Observe that the definitions above are not trivial. For instance, let $\textbf{L}$ be the variety of bounded lattices. It is known (see \cite{V1999} and \cite{FH1970}) that $\textbf{L}$ is a variety with $BFC$. If $L=\textbf{2}\times \textbf{2}$ (with $\textbf{2}$ the chain of two elements) and $M=\{0,1,a,b,c\}$, with $\{a,b,c\}$ not comparables, it easily follows that $L$ is subalgebra of $M$, but $L$ is directly decomposable while $M$ is not. So $\textbf{L}$ is a variety which has not SC nor CSC.
\endrem

Let $\textbf{V}$ be a variety with BFC. We write $\mathcal{V}$ to denote the algebraic category associated to $\textbf{V}$.

\lem\label{preservation imply prods stable by pushouts}
Let $\textbf{V}$ be a variety with BFC. If $\textbf{V}$ has RexDFC and CSC then, in $\mathcal{V}$ the products are stable by pushouts.
\endlem
\pf
Let $A,B\in \textbf{B}$ and $f:A\rightarrow B$ be an homomorphism. If $A\cong A_{1}\times A_{2}$, let us consider de diagram: 

\begin{displaymath}
\xymatrix{
A_{1} \ar[d] & \ar[l] \ar[r] \ar[d] A & A_{2} \ar[d]
\\
P_{1}  & \ar[l] \ar[r]  B & P_{2} 
}
\end{displaymath}
\noindent
Where $P_{1}$ and $P_{2}$ are the pushouts from the left and the right squares, respectively. If $i$ denotes the isomorphism between $A$ and $A_{1}\times A_{2}$, then $A_{j}\cong A/Ker(\pi_{j}i)$ (with $j=1,2$). Since $Ker(\pi_{1}i) \diamond Ker(\pi_{2}i)$ in $Con(A)$, if $\vec{e}_{j}$ denotes the central element corresponding to $Ker(\pi_{j}i)$, thus from Lemmas \ref{Existential implies Positive} and \ref{Positive implies principal congruence} we have that $Ker(\pi_{1}i)=\theta^{A}(\vec{1},\vec{e}_{2})$ and $Ker(\pi_{2}i)=\theta^{A}(\vec{1},\vec{e}_{1})$. From, item 1. of Corollary \ref{corollary diagram pushout}, we get that $P_{1}\cong B/\theta^{B}(\vec{1},f(\vec{e}_{2}))$ and $P_{2}\cong B/\theta^{B}(\vec{1},f(\vec{e}_{1}))$. The universal property of pushouts implies that $B\rightarrow P_{1}$ coincides with $B\rightarrow \theta^{B}(\vec{1},f(\vec{e}_{2}))$ and $B\rightarrow P_{2}$ with $B\rightarrow \theta^{B}(\vec{1},f(\vec{e}_{1}))$. Since $f$ preseserves pairs of complementary central elements by assumption, we can conclude that $B\cong B/\theta^{B}(\vec{1},f(\vec{e}_{2})) \times B/\theta^{B}(\vec{1},f(\vec{e}_{1}))\cong P_{1}\times P_{2}$.
\epf

\lem\label{reciprocal preservation imply prods stable by pushouts}
Let $\textbf{V}$ be a variety with BFC, $A,B\in \textbf{V}$ and $f:A\rightarrow B$ an homomorphism that preserves central elements. If $\textbf{V}$ has RexDFC and in $\mathcal{V}$  the binary products are stable by pushouts along $f$, thus $f$ preserves pairs of complementary central elements.
\endlem
\pf
Let $A\in \mathcal{V}$ and $\vec{e}$ a central element of $A$. If $\vec{g}$ denotes the complementary central element of $\vec{e}$ we get that $\theta^{A}_{\vec{0},\vec{e}}=\theta^{A}(\vec{1},\vec{g})$, so by Lemmas \ref{Existential implies Positive} and \ref{Positive implies principal congruence}, we get that $A\cong A/\theta^{A}(\vec{1},\vec{g})\times A/\theta^{A}(\vec{1},\vec{e})$. Let us, consider the diagram
\begin{displaymath}
\xymatrix{
A/\theta^{A}(\vec{1},\vec{g}) \ar[d] & \ar[l] \ar[r] \ar[d]^-{f} A & A/\theta^{A}(\vec{1},\vec{e}) \ar[d]
\\
B/\theta^{B}(\vec{1},f(\vec{g})) & \ar[l] \ar[r] B & B/\theta^{B}(\vec{1},f(\vec{e}))
}
\end{displaymath}

By Corollary \ref{corollary diagram pushout}, both squares are pushouts, so, since the binary products are stable by pushouts along $f$ by assumption, the span $B/\theta^{B}(\vec{1},f(\vec{g}))\leftarrow B\rightarrow B/\theta^{B}(\vec{1},f(\vec{e}))$ is a product. This fact implies directly that $\theta^{B}(\vec{1},f(\vec{g}))\diamond \theta^{B}(\vec{1},f(\vec{e}))$ in $Con(B)$. Since $f$ preserves central elements by hypothesis, both $f(\vec{e})$ and $f(\vec{g})$ are central elements of $B$ so, we conclude that $f(\vec{e})\diamond_{B} f(\vec{g})$. 
\epf

\lem\label{0 and 1 imply products are codisjoint}
If $\textbf{V}$ is a variety with $\vec{0}$ and $\vec{1}$, then, in $\mathcal{V}$ the pushout of the projections of binary products is the terminal object. 
\endlem
\pf
For every pair of $A,B\in \mathcal{V}$ the pushout of the projections $A\leftarrow A\times B \rightarrow B$ belongs to $\mathcal{V}$. It is clear that the projections send $[\vec{0},\vec{1}]\in A\times B$ to $\vec{0}$ in $A$ and to $\vec{1}$ in $B$, so $\vec{0}=\vec{1}$ in the pushout. Since $\textbf{V}$ is a variety with $\vec{0}$ and $\vec{1}$, it follows that the pushout must be the terminal object.
\epf

\proposition\label{scc is equivalent to coextensivity}
Let $\textbf{V}$ be a variety with BFC. If $\textbf{V}$ has RexDFC with SC, the following are equivalent:
\begin{enumerate}
\item $\textbf{V}$ has CSC.
\item $\mathcal{V}$ is coextensive.
\end{enumerate}
\endproposition
\pf
Since $\textbf{V}$ has BFC, thus is a variety with $\vec{0}$ and $\vec{1}$. So, by Lemma \ref{0 and 1 imply products are codisjoint}, we get that the pushout of the projections of binary products is the terminal object. Let us assume that $\textbf{V}$ has SCC. From the Lema \ref{preservation imply prods stable by pushouts}, the products are stable by pushouts. Hence, by the dual of the Proposition \ref{Coextensivity Proposition}, $\mathcal{V}$ is coextensive. The  reciprocal follows from Lemma \ref{reciprocal preservation imply prods stable by pushouts}.
\epf

\section{An axiomatization for connected models}\label{Axiomatization of connected models}

In this section we prove that a variety with BFC has RexDFC (LexDFC) if and only if the factor congruence $\theta_{\vec{1},\vec{e}}$ associated to a central element $\vec{e}$, coincides with the principal congruence that identifies $\vec{1}$ with $\vec{e}$ ($\vec{0}$ with $\vec{e}$). This fact will allows us to prove that the theory of connected models for varieties with RexDFC (LexDFC) is definable by a finite set of first order formulas.
\\

We will use the following (Gr\"atzer) version of Maltsev's key observation on principal congruences.

\lem\label{Gratzer Malsev Lemma} Let $A$ be an algebra and $a, b \in \textbf{A}$, $\vec{c}$, $\vec{d}$ $\in A^{n}$. Then $(a,b)\in \theta^{A}(\vec{c},\vec{d})$ if and only if there exist $(n+m)$-ary terms $t_{1}(\vec{x},\vec{u})$,...,$t_{k}(\vec{x},\vec{u})$ with $k$ odd and $\vec{\lambda}\in A^{m}$ such that: 

\begin{center}
\begin{tabular}{cc}
$a=t_{1}(\vec{c},\vec{\lambda})$ & $b=t_{k}(\vec{d},\vec{\lambda})$
\\
$t_{i}(\vec{c},\vec{\lambda})=t_{i+1}(\vec{c},\vec{\lambda})$, $i$ even, & $t_{i}(\vec{d},\vec{\lambda})=t_{i+1}(\vec{d},\vec{\lambda}),$ $i$ odd.
\end{tabular}
\end{center}
\endlem

We recall that a \emph{principal congruence formula} is a formula $\pi(x,y,\vec{u},\vec{v})$ of the form
\[\exists_{\vec{w}}(x\approx t_{1}(\vec{u},\vec{w})\wedge \bigwedge_{i\in E_{k}} (t_{i}(\vec{u},\vec{w})\approx t_{i+1}(\vec{u},\vec{w}))\wedge \bigwedge_{i\in O_{k}}(t_{i}(\vec{v},\vec{w})\approx t_{i+1}(\vec{v},\vec{w}))\wedge t_{k}(\vec{v},\vec{w})\approx y)) \]

where $k$ is odd and $t_{i}$ are terms of type $\tau$. This fact allows us to restate the latter Lemma as 

\lem\label{Malsev restated}
Let $A$ be an algebra, $a, b \in \textbf{A}$, $\vec{c}$, $\vec{d}$ $\in A^{n}$. Then $(a,b)\in \theta^{A}(\vec{c},\vec{d})$ if and only if there exists a principal congruence formula $\pi$, such that $A\models \pi(a,b,\vec{c},\vec{d})$.
\endlem

\lem\label{Definability by principal congruences}
Let $\textbf{V}$ be a variety with DFC, $A\in \textbf{V}$ and $\vec{e}\in Z(A)$: 
\begin{enumerate}
\item If $\theta^{A}_{\vec{1},\vec{e}}=\theta^{A}(\vec{1},\vec{e})$ then $\theta_{\vec{1},\vec{e}}$ is definible by a formula of the form $\exists \bigwedge p\approx q$. 
\item If $\theta^{A}_{\vec{0},\vec{e}}=\theta^{A}(\vec{0},\vec{e})$ then $\theta_{\vec{0},\vec{e}}$ is definible by a formula of the form $\exists \bigwedge p\approx q$.
\end{enumerate}

\endlem
\pf
We only prove $1.$ because the proof of $2.$ is essentially the same. Let us write $P=\textbf{F}(x,y)\times \textbf{F}(y)$, where $\textbf{F}(x,y)$ and $\textbf{F}(y)$ are the free algebras generated by $\{x,y\}$ and $\{y\}$, respectively. By hyphotesis, $Ker(\pi_{2})=\theta^{\bf{P}}_{[\vec{1},\vec{1}],[\vec{0},\vec{1}]}=\theta^{P}([\vec{1},\vec{1}],[\vec{0},\vec{1}])$. Since the pair $((x,y),(y,y))\in Ker(\pi_{2})$, from Lema \ref{Gratzer Malsev Lemma}, there exist $(n+m)$-ary terms $t_{1}(\vec{x},\vec{u})$,...,$t_{k}(\vec{x},\vec{u})$ with $k$ odd and $\vec{u}\in P^{m}$ such that:

\begin{equation}
\begin{array}{cc}
(x,y)=t^{P}_{1}[[\vec{1},\vec{1}],\vec{u}]& (y,y)=t^{P}_{k}[[\vec{0},\vec{1}],\vec{u}]
\\
t^{P}_{i}[[\vec{1},\vec{1}],\vec{u}]=t^{P}_{i+1}[[\vec{1},\vec{1}],\vec{u}],\; i\in E_{k}, & t^{P}_{i}[[\vec{0},\vec{1}],\vec{u}]=t^{P}_{i+1}[[\vec{0},\vec{1}],\vec{u}],\; i\in O_{k}.\label{1}
\end{array}
\end{equation}

where $E_{k}$ and $O_{k}$ refer to the even and odd naturals less or equal to $k$, respectively.
\\

Since $\vec{u}\in P$, there are $\vec{P}(x,y)\in F(x,y)$ and $\vec{Q}(y)\in F(x,y)$, such that $\vec{u}=[\vec{P}, \vec{Q}]$. Recall that $t^{P}_{i}[[\vec{R},\vec{S}],[\vec{P}, \vec{Q}]]=(t^{F(x,y)}_{i}[\vec{R},\vec{P}],t^{F(y)}_{i}[\vec{S},\vec{Q}])$, for $1\leq i \leq k$ and $[\vec{R},\vec{S}]\in P$, thus, from equation (\ref{1}), we obtain that there exist $(n+m)$-ary terms $t_{1}(\vec{x},\vec{u})$,...,$t_{k}(\vec{x},\vec{u})$ with $k$ odd, $\vec{P}(x,y)\in F(x,y)$ and $\vec{Q}(y)\in F(y)$, such that:  

{\small
\begin{displaymath}
y=t^{F(y)}_{i}[\vec{0},\vec{Q}(y)], \textrm{for every}\; 1\leq i\leq k 
\end{displaymath}
}
and

{\small
\begin{displaymath}
\begin{array}{cc}
x=t^{F(x,y)}_{1}[\vec{1},\vec{P}(x,y)] & y=t^{F(x,y)}_{k}[\vec{0},\vec{P}(x,y)]
\\
t^{F(x,y)}_{i}[\vec{1},\vec{P}(x,y)]=t^{F(x,y)}_{i+1}[\vec{1},\vec{P}(x,y)],\; i\in E_{k} & t^{F(x,y)}_{i}[\vec{0},\vec{P}(x,y)]=t^{F(x,y)}_{i+1}[\vec{0},\vec{P}(x,y)],\; i\in O_{k}
\end{array}
\end{displaymath}
}  
  
Let $\varphi(x,y,\vec{z})=\pi(x,y,\vec{1},\vec{z})$. In order to to check that $\varphi$ defines $\theta^{A}_{\vec{1},\vec{e}}$ in terms of $\vec{e}$ let us assume $A,B\in \mathcal{V}$ and $(a,b),(c,d)\in A\times B$. Since the free algebra functor $\textbf{F}:\Set \rightarrow \mathcal{V}$ is left adjoint to the forgetful functor, in the case of $b=d$, the assingments $\alpha_{A}:\{x,y\}\rightarrow A$ and $\alpha_{B}:\{y\}\rightarrow B$, defined by $\alpha_{A}(x)=a$, $\alpha_{A}(y)=c$ and $\alpha_{B}(y)=b$ generate a unique pair of homomorphisms $\beta_{A}:\textbf{F}(x,y)\rightarrow A$ and $\beta_{B}:\textbf{F}(y)\rightarrow B$ extending $\alpha_{A}$ and $\alpha_{B}$, respectively. Therefore, since $P\models \varphi((x,y),(y,y), [\vec{0},\vec{1}])$, by applying $g=\beta_{A}\times \beta_{B}$ in (\ref{1}), we obtain as result that $A\times B\models \varphi((a,b),(c,b),[\vec{0},\vec{1}])$. On the other hand, if $A\times B\models \varphi((a,b),(c,d),[\vec{0},\vec{1}])$, then there exist $[\vec{\varepsilon},\vec{\delta}]\in A\times B$, such that 

{\small
\begin{equation}\label{2}
\begin{array}{cc}
(a,b)=t^{A\times B}_{1}[[\vec{1},\vec{1}],[\vec{\varepsilon},\vec{\delta}]] & (c,d)=t^{A\times B}_{k}[[\vec{0},\vec{1}],[\vec{\varepsilon},\vec{\delta}]]
\\
t^{A\times B}_{i}[[\vec{1},\vec{1}],[\vec{\varepsilon},\vec{\delta}]]=t^{A\times B}_{i+1}[[\vec{1},\vec{1}],[\vec{\varepsilon},\vec{\delta}]],\;i\in E_{k}, & t^{A\times B}_{i}[[\vec{0},\vec{1}],[\vec{\varepsilon},\vec{\delta}]]=t^{A\times B}_{i+1}[[\vec{0},\vec{1}],[\vec{\varepsilon},\vec{\delta}]],\;i\in O_{k}.
\end{array}
\end{equation}
}

So, since $t^{A\times B}_{i}[[\vec{j},\vec{r}],[\vec{\varepsilon},\vec{\delta}]]=(t^{A}_{i}[\vec{j},\vec{\varepsilon}],t^{B}_{i}[\vec{r},\vec{\delta}])$, for every $[\vec{j},\vec{r}]\in A\times B$ and $1\leq i\leq k$, from (\ref{2}), we conclude that

\begin{displaymath}
\begin{array}{cc}
b=t^{B}_{1}[\vec{1},\vec{\varepsilon}] & d=t^{B}_{k}[\vec{1},\vec{\varepsilon}]
\\
t^{B}_{i}[\vec{1},\vec{\varepsilon}]=t^{B}_{i+1}[\vec{1},\vec{\varepsilon}],\; i\in E_{k} & t^{B}_{i}[\vec{1},\vec{\varepsilon}]=t^{B}_{i+1}[\vec{1},\vec{\varepsilon}],\; i\in O_{k}.
\end{array}
\end{displaymath}

Which by Lemma \ref{Gratzer Malsev Lemma} means that $(b,d)\in \theta^{B}(\vec{1},\vec{1})$. Since $\theta^{B}(\vec{1},\vec{1})=\theta^{B}_{\vec{1},\vec{1}}$ by assumption and $\theta^{B}(\vec{1},\vec{1})=\Delta^{B}$, we get that $b=d$. This concludes the proof.

\epf

\corollary\label{Characterization of varieties with RexDFC}
Let $\V$ be a variety with DFC. The following are equivalent:
\begin{enumerate}
\item $\V$ has RexDFC if and only if, for every $A\in \textbf{V}$ and $\vec{e}\in Z(A)$, $\theta^{A}_{\vec{1},\vec{e}}=\theta^{A}(\vec{1},\vec{e})$.
\item $\V$ has LexDFC if and only if, for every $A\in \textbf{V}$ and $\vec{e}\in Z(A)$, $\theta^{A}_{\vec{0},\vec{e}}=\theta^{A}(\vec{0},\vec{e})$.
\end{enumerate}
\endcorollary
\pf
In each item, the first implication follows from Lemmas \ref{Existential implies Positive}, \ref{Positive implies principal congruence} and the last one is a consequence of Lemma \ref{Definability by principal congruences}.
\epf

We say that a set of formulas $\Sigma(\vec{z},\vec{u})$ \emph{defines the property} $\vec{e}\diamond_{A} \vec{f}$ \emph{in} $\textbf{V}$ if for every $A\in \textbf{V}$ and $\vec{e},\vec{f}\in A^{n}$ it follows that $\vec{e}\diamond_{A} \vec{f}$ if and only if $A\models \sigma[\vec{e},\vec{f}]$, for every $\sigma\in \Sigma$.
\\

Let $\vec{e},\vec{f}\in A^{n}$ and $\varphi$ be the formula used in the proof of Lemma \ref{Definability by principal congruences}. We consider the following formulas:
\\

\noindent
$\tau_{r}(\vec{z},\vec{u})=(\forall_{x})(\varphi(x,x,\vec{z}))$
\\
$\tau_{s}(\vec{z},\vec{u})=(\forall_{x,y})(\varphi(x,y,\vec{z})\rightarrow \varphi(y,x,\vec{z}))$
\\
$\tau_{t}(\vec{z},\vec{u})=(\forall_{x,y,v})(\varphi(x,v,\vec{z})\wedge \varphi(v,y,\vec{z})\rightarrow \varphi(x,y,\vec{z}))$
\\
$\tau_{i}(\vec{z},\vec{u})=(\forall_{x,y})(\varphi(x,y,\vec{z})\wedge \varphi(x,y,\vec{u})\rightarrow x\approx y)$
\\
$\tau_{p}(\vec{z},\vec{u})=(\forall_{x,y})(\exists_{v})(\varphi(x,v,\vec{z})\wedge \varphi(v,y,\vec{u}))$
\\
$\tau_{k}(\vec{z},\vec{u})=\bigwedge_{1\leq  j\leq  n}\varphi(1_{j},z_{j},\vec{z})\wedge \bigwedge_{1\leq  j\leq  n}\varphi(0_{j},u_{j},\vec{z})$
\\

And for every function symbol $f$ in the languange of $A$: 
\\

\noindent
$\tau_{f}(\vec{z},\vec{u})=(\forall_{l_{1},...,l_{m},v_{1},...,v_{m}})(\bigwedge_{1\leq \alpha\leq m}\varphi(l_{\alpha},v_{\alpha},\vec{z})\rightarrow \varphi(f(l_{1},...,l_{m}),f(v_{1},...,v_{m}),\vec{z}))$
\\

If we call $E_{0}=\{\tau_{\beta}(\vec{z},\vec{u})\mid \beta\in \{r,s,t,i,p,k\}\}$, $E_{1}=\{\tau_{\beta}(\vec{u},\vec{z})\mid \beta\in \{r,s,t,i,p,k\}\}$, $C= \{\tau_{f}(\vec{z},\vec{u})\mid f\in \tau\}$, where $\tau$ is the type of $A$, let $\Sigma(\vec{z},\vec{u})=E_{0}\cup E_{1}\cup C$.

\lem\label{Definability of beeing factor}
Let $\textbf{V}$ a variety with BFC. 
\begin{enumerate}
\item If $\textbf{V}$ has RexDFC, there exists a set of formulas $\Sigma(\vec{z},\vec{u})$ defining the property $\vec{e}\diamond_{A} \vec{f}$ in $\textbf{V}$.
\item If $\textbf{V}$ has LexDFC, there exists a set of formulas $\Sigma(\vec{z},\vec{u})$ defining the property $\vec{e}\diamond_{A} \vec{f}$ in $\textbf{V}$.
\end{enumerate}
\endlem
\pf
We prove $1.$ Let $A\in \textbf{V}$, $\vec{e},\vec{f}\in A^{n}$ and suppose that $\textbf{V}$ has LexDFC. We define the following relations in  $A$:
\begin{center}
\begin{tabular}{cc}
$L_{\vec{e}}=\{(a,b)\in A\times A\mid A\models \varphi[\vec{e},a,b]\}$ & $L_{\vec{f}}=\{(a,b)\in A\times A\mid A\models \varphi[\vec{f},a,b]\}$
\end{tabular}
\end{center}

Let $\vec{e},\vec{f}\in A^{n}$. Observe that formulas $\tau_{r}$, $\tau_{s}$ and $\tau_{t}$ say that $L_{\vec{e}}$ is an equivalence relation on $A$. The set $\{\tau_{f}\mid \textrm{f is a symbol of function in the lenguage of A}\}$ says is that $L_{\vec{e}}$ is a congruence. The formula $\tau_{i}$ says that $L_{\vec{e}} \cap L_{\vec{f}}=\Delta^{A}$ and the formula $\tau_{i}$ says that $L_{\vec{e}} \circ L_{\vec{f}}=\nabla^{A}$. Finally, the formula $\tau_{k}$ says that $[\vec{1},\vec{e}]\in L_{\vec{e}}$ and $[\vec{0},\vec{f}]\in L_{\vec{e}}$. 
\\

It is clear that if $A\models \sigma[\vec{e},\vec{f}]$ for every $\sigma \in \Sigma$, then $L_{\vec{e}}$ and $L_{\vec{f}}$ are factor congruences of $A$ such that $\vec{e}\equiv \vec{0} (L_{\vec{e}})$, $\vec{f}\equiv \vec{1} (L_{\vec{e}})$, $\vec{f}\equiv \vec{0} (L_{\vec{f}})$, $\vec{e}\equiv \vec{1} (L_{\vec{f}})$. That is, $L_{\vec{e}}=\theta^{A}_{\vec{0},\vec{e}}$ and $L_{\vec{f}}=\theta^{A}_{\vec{0},\vec{f}}$. Hence, $\vec{e},\vec{f}\in Z(A)$ and $\vec{e}\diamond_{A}\vec{f}$. On the other hand, if $\vec{e}\diamond_{A}\vec{f}$, from Lemmas \ref{Existential implies Positive}, \ref{Positive implies principal congruence} and \ref{Definability by principal congruences}, we get that $L_{\vec{e}}=\theta^{A}_{\vec{0},\vec{e}}=\theta^{A}(\vec{0},\vec{e})$ and $L_{\vec{f}}=\theta^{A}_{\vec{0},\vec{f}}=\theta^{A}(\vec{0},\vec{f})$, so $A\models \sigma[\vec{e},\vec{f}]$ for every $\sigma \in \Sigma$. The proof of $2.$ is similar.
\epf

Again, let $\V$ be a variety with BFC. We write $\V_{C}$ to denote the class of connected (directly idecomposable) algebras of $\V$. If $A\in \V_{C}$, then we also say that $A$ is a $\V$-connected algebra. 

\corollary\label{Theory of connected models}
If $\textbf{V}$ has RexDFC (or LexDFC), the class $\textbf{V}_{C}$ is axiomatizable by a set of first order formulas.
\endcorollary
\pf
Suppose that $\textbf{V}$ has RexDFC. Consider the set $\Sigma$ from Proposition \ref{Definability of beeing factor}. It is immidiate that $A\in \textbf{V}_{DI}$ if and only if in $A$ the following axioms hold
\begin{center}
$\vec{0}\neq \vec{1}$ and $\forall_{\vec{e},\vec{f}}\bigwedge \Sigma(\vec{e},\vec{f})\rightarrow ((\vec{e}=\vec{0} \wedge \vec{f}=\vec{1})\vee (\vec{e}=\vec{1} \wedge \vec{f}=\vec{0})).$ 
\end{center} 
\epf

\section{Connected models in a Topos}

Let $\V$ be a variety and $\Sigma(\vec{x},\vec{y})$ the set of formulas of item $1.$ in Lemma \ref{Definability of beeing factor}. We call $\Vmod$ to the theory given by the equations holding in $\V$ and the axiom 
\begin{equation}\label{Axiom 0 and 1}
\top \vdash \bigwedge \Sigma (\vec{1},\vec{0})
\end{equation}

For a given topos $\mathcal{E}$, let $\Vmod(\mathcal{E})$ be the category of internal models in $\mathcal{E}$ respect to $\Vmod$. Observe that in $\Set$, axiom (\ref{Axiom 0 and 1}) is equivalent to say that $\V$ is a variety with $\vec{0}$ and $\vec{1}$.
\\

With the aim of understand what is a variety with RexDFC in a topos $\mathcal{E}$, observe that the proof of Lemma \ref{Definability by principal congruences} suggest that we can get a weaker condition to make a variety with $\vec{0}$ and $\vec{1}$ be variety with BFC. That is: \emph{Let $\V$ be a variety with $\vec{0}$ and $\vec{1}$; and let $A,B\in \V$. Consider the projection $\pi_{B}:A\times B\rightarrow B$. If $Ker(\pi_{B})=\theta^{A\times B}([\vec{1},\vec{1}],[\vec{0},\vec{1}])$ for every $A,B\in \V$, then $\V$ has BFC. Furthermore, $\V$ has RexDFC.} The proof of this fact uses the same arguments of the one given for Lemma \ref{Definability by principal congruences} in order to obtain an existential formula defining $Ker(\pi_{B})$. This entails that $\V$ has (DFC) and consecuently (see the Introduction) $\V$ has (BFC).
\\

So, let $A$ be in $\Vmod(E)$ and for every $1\leq i\leq n$ consider the following composites:
\begin{displaymath}
\xymatrix{
1\ar[r]^-{\vec{0}} \ar@/_1pc/[rr]_-{0_{i}} & A^{n} \ar[r]^-{\pi_{i}} & A & & 1\ar[r]^-{\vec{1}} \ar@/_1pc/[rr]_-{1_{i}} & A^{n} \ar[r]^-{\pi_{i}} & A
}
\end{displaymath}

Thus, for every $i$ and $B$ in $\Vmod(E)$ we obtain a morphism

\begin{displaymath}
\xymatrix{
1 \ar[r]^-{f_{i}} &  (A\times B) \times (A\times B)
}
\end{displaymath}

where $f_{i}=\langle \langle 1_{i},1_{i}\rangle, \langle 0_{i},1_{i}\rangle \rangle$. Consider the projection $\pi_{B}:A\times B \rightarrow B$ and let $b:Ker(\pi_{B})\rightarrow (A\times B)^{2}$ be the morphism induced by the span $A\times B \leftarrow Ker(\pi_{B})\rightarrow A\times B$. It easily follows that every $f_{i}$ factors through $b$.

\definition\label{RexDFC in a topos}
A category of internal $\Vmod$-models in $\mathcal{E}$ has the {\Rex} (KDP) if for every pair of objects $A,B$ in $\Vmod(E)$, $Ker(\pi_{B})$ is the least subobject of $(A\times B)^{2}$ in $\Vmod(\mathcal{E})$ through which the collection $\{f_{i}\mid 1\leq i\leq n\}$ factors. I,e:
\begin{displaymath}
\xymatrix{
1 \ar[rr]^-{f_{i}} \ar[dr]_-{l_{i}} \ar@/_1pc/[ddr]_-{a_{i}} & & (A\times B)^{2}
\\
 & C \ar[ur]_-{m} & 
\\
 & Ker(\pi_{B})\ar@/_1pc/[uur]_-{b} \ar@{-->}[u]_-{k} & 
}
\end{displaymath}  
If $m:C\rightarrow (A\times B)^{2}$ is a subobject such that $ml_{i}=f_{i}$ for every $1\leq i\leq n$, then there exists a morphims $k:Ker(\pi_{B})\rightarrow C$ (necessarily unique), such that $mk=b$ and $ka_{i}=l_{i}$.
\enddefinition


Inspired in Corollary \ref{Theory of connected models} we introduce the following definition

\definition\label{Connected VModels in a topos}
Let $\Vmod(E)$ a category of internal $\Vmod$-models in $\mathcal{E}$ with KDP. An internal $\Vmod$-model $A$ is connected if the following sequents hold 
\begin{center}
\begin{itemize}
\item[(C1)] $\vec{0}=\vec{1} \vdash \perp$
\item[(C2)] $\bigwedge \Sigma (\vec{x},\vec{y}) \vdash_{\vec{x},\vec{y}} (\vec{x}=\vec{0} \wedge \vec{y}=\vec{1}) \vee (\vec{x}=\vec{1} \wedge \vec{y}=\vec{0}) $
\end{itemize}
\end{center}
in the internal logic of $\mathcal{E}$.
\enddefinition

In the following, we write $\CV$ for the theory of internal connected $\Vmod$-models in $\mathcal{E}$. 
\\

Suppose $A$ is in $\CV$. Observe that, from axiom (\ref{Axiom 0 and 1}), there exists a morphism $g$ such that the diagram below

\begin{displaymath}
\xymatrix{
1 \ar[r]^-{\langle \vec{1},\vec{0}\rangle} \ar[d]  & A^{n}\times A^{n} 
\\
Im(\langle \vec{1},\vec{0}\rangle) \ar[r]_-{g}  & [\bigwedge \Sigma (\vec{x},\vec{y})] \ar[u]
}
\end{displaymath}

commutes. Since $Im(\langle \vec{1},\vec{0}\rangle)\cong [\vec{x}=\vec{1} \wedge \vec{y}=\vec{0}]$, we get that $[\vec{x}=\vec{1} \wedge \vec{y}=\vec{0}]\leq [\bigwedge \Sigma (\vec{x},\vec{y})]$ in $Sub(A^{n}\times A^{n})$. Thus, by completeness (c.f. D1.4.11 in \cite{J2002}) we get that the sequent
$(\vec{x}=\vec{1} \wedge \vec{y}=\vec{0})\vdash_{\vec{x},\vec{y}} [\bigwedge \Sigma (\vec{x},\vec{y})]$ holds in the internal logic of $\mathcal{E}$. Moreover, since \[Im(\langle \vec{1},\vec{0}\rangle)\cong Im(\langle \vec{0},\vec{1}\rangle)\cong 1\cong [\vec{x}=\vec{0} \wedge \vec{y}=\vec{1}]\]
from axiom (\ref{Axiom 0 and 1}), $Im(\langle \vec{0},\vec{1}\rangle)\leq [\Sigma (\vec{x},\vec{y})]$ in $Sub(A^{n}\times A^{n})$, so the sequent $\top \vdash \bigwedge \Sigma (\vec{0},\vec{1})$ holds in the internal logic of $\mathcal{E}$. By proceeding as before, we can deduce that the sequent $(\vec{x}=\vec{0} \wedge \vec{y}=\vec{1})\vdash_{\vec{x},\vec{y}} [\bigwedge \Sigma (\vec{x},\vec{y})]$ also holds in the internal logic of $\mathcal{E}$. We have proved the following

\lem\label{Non obvious sequent}
In $\CV$ the sequent $(\vec{x}=\vec{1} \wedge \vec{y}=\vec{0})\vee (\vec{x}=\vec{0} \wedge \vec{y}=\vec{1}) \vdash_{\vec{x},\vec{y}} [\bigwedge \Sigma (\vec{x},\vec{y})]$ holds.
\endlem

Let us consider the points $\langle \vec{0},\vec{1}\rangle: 1\rightarrow A^{n}\times A^{n}$ and $\langle \vec{1},\vec{0}\rangle: 1\rightarrow A^{n}\times A^{n}$. From axiom (\ref{Axiom 0 and 1}) it follows that, for every $\sigma \in \Sigma$, there exist a morphism $l_{\sigma}:1+1\rightarrow [\sigma(\vec{x},\vec{y})]$, such that the diagram below
\begin{displaymath}
\xymatrix{
1+1 \ar[r]^-{l_{\sigma}}  \ar[dr]_-{[\langle \vec{0},\vec{1}\rangle,\langle \vec{1},\vec{0}\rangle]} & [\sigma(\vec{x},\vec{y})] \ar[d]
\\
 & A^{n}\times A^{n}
}
\end{displaymath}
commutes, so $1+1\leq [\sigma(\vec{x},\vec{y})]$ for every $\sigma \in \Sigma$.  Hence, $1+1\leq [\bigwedge \Sigma(\vec{x},\vec{y})]$ in $Sub(A^{n}\times A^{n})$. Let us call $\alpha: 1+1\rightarrow [\Sigma(\vec{x},\vec{y})]$ to the morphism that arise from the factorization of $1+1 \rightarrow A^{n}\times A^{n}$ along $[\Sigma(\vec{x},\vec{y})] \rightarrow A^{n}\times A^{n}$.
\\

Since $[(\vec{x} = \vec{1} \wedge \vec{y} = \vec{0}) \vee (\vec{x} = \vec{0} \wedge \vec{y} = \vec{1})] \cong 1 + 1$, as result of the latter discution, we obtain a characterization for the internal connected $\Vmod$-models in $\mathcal{E}$.

\lem\label{Characterization of internal models in a topos}
Let $\Vmod(E)$ a category of internal $\Vmod$-models in $\mathcal{E}$ with KDP. An internal $\Vmod$-model $A$ in $\mathcal{E}$ is connected if and only if the diagram below 
\begin{displaymath}
\xymatrix@1{
0 \ar[r]^-{!}  & 1 \ar@<0.5ex>[r]^-{\vec{1}} \ar@<-0.5ex>[r]_-{\vec{0}}  & A^{n} 
}
\end{displaymath}
is an equalizer in $\mathcal{E}$, and the morphism $\alpha: 1+1\rightarrow [\Sigma(\vec{x},\vec{y})]$ is an iso.
\endlem
\pf
Since $1+1\leq [\Sigma(\vec{x},\vec{y})]$ in $Sub(A^{n}\times A^{n})$, the result follows from apply Lemma \ref{Non obvious sequent} and the interpretations of the axioms (C1) and (C2) of Definition \ref{Connected VModels in a topos} in the internal logic of $\mathcal{E}$.
\epf

\section{Connected models in Coherent topoi}

It is known that every distributive lattice $D$ can be treaten as a \emph{coherent category} (A1.4 in \cite{J2002}). Its \emph{coherent coverage} (A2.1.11(b) in \cite{J2002}) is the function that sends each $d \in D$ to the set of finite families $\{d_{i} \leq d \mid i \in I\}$ such that $\bigvee_{i\in I} d_{i} = d$. As usual, the resulting topos of sheaves will be denoted by $\Shv(D)$. Binary covers $a \vee b = d$ of $d \in D$ will play an important role because in order to check that a presheaf $P : D^{op} \rightarrow \Set$ is a sheaf, it is enough to check the sheaf condition for binary covers.

Recall that every variety $\V$, is an algebrabic category over $\Set$. Thus a $\V$-model in $\Shv(D)$ is a functor $D^{op}\rightarrow \V$ such that the composite presheaf $D^{op}\rightarrow \V \rightarrow \Set$ is a sheaf.

The aim of this section is to characterize internal connected $\Vmod$-models in $\Shv(D)$. Since the set of formulas $\Sigma$ of Lemma \ref{Definability by principal congruences} is not composed by equations, the theory $\Vmod$ is not algebraic. So, we need to understand first what an internal $\Vmod$-model in $\Shv(D)$ is. Observe that, from axiom (\ref{Axiom 0 and 1}), it follows that a $\V$-model $A$ in $\Shv(D)$ is an internal $\Vmod$-model in $\Shv(D)$ if and only if axiom (\ref{Axiom 0 and 1}) holds in $A(d)$, for every $d\in D$. Hence an internal $\Vmod$-model in $\Shv(D)$ is just a sheaf $A$ such that $A(d)$ is an algebra of $\V$ with $\vec{0}$ and $\vec{1}$, for every $d\in D$.

\lem\label{RexDFC makes V(PSh(D)) having KDP}
Let $\V$ be a variety with $\vec{0}$ and $\vec{1}$ in $\Set$. If $\V$ has RexDFC then $\Vmod(\widehat{D})$ has KDP.
\endlem
\pf
Since $\V$ is a variety with $\vec{0}$ and $\vec{1}$ then the axiom (\ref{Axiom 0 and 1}) holds. Let $F,G$ in $\widehat{D}$. Since limits in $\widehat{D}$ are calculated pointwise, then if $\pi_{G}:F\times G\rightarrow G$, it follows that $Ker(\pi_{G})(d)=Ker(\pi_{G(d)})$, for every $d\in D$. Thus, since $\V$ has RexDFC by hypothesis, 
\[Ker(\pi_{G})(d)=\theta^{F(d)\times G(d)}([\vec{1}_{F(d)},\vec{1}_{G(d)}],[\vec{0}_{F(d)},\vec{1}_{G(d)}])=\theta^{F\times G}([\vec{1},\vec{1}],[\vec{0},\vec{1}])(d)\]
That is, $Ker(\pi_{G})$ is the least subobject of $(F\times G)^{2}$ in $\Vmod(\widehat{D})$ through which the collection $\{f_{i}\mid 1\leq i\leq n\}$ factors. This concludes the proof.
\epf

\corollary\label{RexDFC makes V(Sh(D)) having KDP}
Let $\V$ be a variety with $\vec{0}$ and $\vec{1}$ in $\Set$. If $\V$ has RexDFC then $\Vmod(\Shv(D))$ has KDP.
\endcorollary
\pf
Let $F,G$ in $\Shv(D)$, and let $H\rightarrow (F\times G)^{2}$ be a subobject in $\Shv(D)$ such that upper triangle in the diagram below  
\begin{displaymath}
\xymatrix{
1 \ar[rr]^-{f_{i}} \ar[dr]_-{l_{i}} \ar@/_1pc/[ddr]_-{a_{i}} & & (F\times G)^{2}
\\
 & H \ar[ur]_-{m} & 
\\
 & Ker(\pi_{G})\ar@/_1pc/[uur]_-{b} \ar@{-->}[u]_-{k} & 
}
\end{displaymath}  

commutes in $\widehat{D}$, for every $1\leq i\leq n$. From Lemma \ref{RexDFC makes V(PSh(D)) having KDP}, $\Vmod(\widehat{D})$ has KDP, so there exists a unique $k:Ker(\pi_{G})\rightarrow H$, such that the left and the right triangles in the diagram above, commutes. Let $\textbf{a}:\widehat{D}\rightarrow \Shv(D)$ be the sheafification respect to the coherent site. Since $\textbf{a}$ preserves finite limits, then $\textbf{a}(Ker(\pi_{G}))$ is the kernel of $\pi_{G}$ in $\Shv(D)$. This concludes the proof.
\epf

It is clear that $\mathbb{C}_{\Vmod}$ neither is an algebraic theory; so to achieve our goal, we will require a  little more effort. To do so, we will use specifically, a suitable description of binary coproducts in $\Shv(D)$ proved in (\cite{CMZ2016}). 

\lem[Binary coproducts in $\Shv(D)$]\label{LemCoprosInShvD} For every $X$, $Y$ in $\Shv(D)$, the coproduct ${X + Y}$ may be defined by 
\[ (X + Y)(d) = \{ (a, b, x, y) \mid a \vee b = d , a\wedge b = \bot, x\in X (a) , y \in Y (b) \} \]
and, for any ${(a, b, x, y)\in (X + Y) (d)}$, 
\[ (X + Y)(c\leq d)(a, b, x, y)  = (a, b, x, y) \cdot c = (a\wedge c, b\wedge c, x\cdot (a\wedge c), y\cdot (b\wedge c)) \]
where ${x\cdot c = X(c\leq d)(x) \in X (c)}$ and ${y\cdot c = Y(c\leq d)(y) \in Y (c)}$.
\endlem

In particular, from Lemma \ref{LemCoprosInShvD}, $(1+1)(d)=\{(a,b)\mid a\vee b=d, a\wedge b=0\}$. That is, $1+1$ is the ``object of partitions" of $D$. 
\\

Recall that the variety of bounded distributive lattices is a variety with DFC, so for a bounded distributive lattice $D$ the subobject $(1+1)(d)$ is isomorphic to the set $Z(\downarrow d)$; i.e, the center of $\downarrow d=\{a\in D\mid a\leq d\}$. On the other hand, if $\V$ is a variety with RexDFC, Corollary \ref{RexDFC makes V(Sh(D)) having KDP}  tells us that $\Vmod(\Shv(D))$ has KDP, so if $A:D^{op}\rightarrow \Set$ is an internal $\Vmod$-model in $\Shv(D)$, then, $A(d)$ has BFC for every $d\in D$. This fact combined with Corollary \ref{Theory of connected models} allows to say that, for every $d\in D$, the interpretation of $[\bigwedge \Sigma(\vec{x},\vec{y})](d)$ brings an isomorphic description of $Z(A(d))$. Finally, from Lemma \ref{Characterization of internal models in a topos}, we get that $\alpha :1+1\rightarrow [\bigwedge \Sigma(\vec{x},\vec{y})]$ is a natural iso, so for every $d\in D$, the map $\alpha_{d}:Z(\downarrow d)\rightarrow Z(A(d))$ is bijective.  

\proposition\label{connected models in ShvD}
Let $\V$ be a variety with RexDFC. An internal $\Vmod$-model $A$ is connected in $\Shv(D)$ if and only if the following conditions hold:
\begin{enumerate}
\item If $A(d)=1$ then $d=0$.
\item For every $c\leq d\in D$, the diagram below
\begin{displaymath}
\xymatrix{
Z(A(d)) \ar[r]^-{\alpha_{d}} \ar[d]_-{k_{d}} & Z(\downarrow d) \ar[d]^-{j_{d}}
\\
Z(A(c)) \ar[r]_-{\alpha_{c}} & Z(\downarrow c)
}
\end{displaymath}
commutes, where, for every $c\in D$, $j_{c}(a)=a \wedge c $, $k_{c}(\vec{e})=\vec{e}\cdot c$ and $\alpha_{c}$ is an isomorphism.  

\end{enumerate}
\endproposition
\pf
A restatement of Lemma \ref{Connected VModels in a topos} in the case of $\Shv(D)$.
\epf

\section{The category of representations}\label{The category of representations}

\definition\label{Category of Representations}
A representation (of a $\textbf{V}$-model) is a pair $(D,X)$, consisting of a distributive lattice $D$ and a $\textbf{V}$-model in $\Shv(D)$ satisfying the equivalent conditions of 	Proposition \ref{connected models in ShvD}. 
\enddefinition

We now define a category $\Rep$ whose objects are representations in the above sense. To describe the arrows in $\Rep$ first recall that any morphism $f : D \rightarrow E$ between distributive lattices is in fact a morphism of sites (in the sense of Theorem VII.10.1 in \cite{MM2012}) when $D$ and $E$ are considered as small categories equipped with the coherent topology. From Theorem VII.10.1 in \emph{loc. cit.} there exists a geometric morphism $f: \Shv(E) \rightarrow \Shv(D)$, whose direct image $f_{\ast}$ is defined as $f_{\ast}(X)=X\circ f$, for every $X\in \Shv(D)$; i.e. $f_{\ast}(X)(d)=X(f(d))$, for every $d\in D$. 
\\

We define now the maps in $\Rep$. For representations $(D, X)$ and $(E, Y )$, an arrow $(D, X) \rightarrow (E, Y )$ in $\Rep$ is a pair $(f, \varphi)$ with $f : D \rightarrow E$ a morphism in $\dLat(\Set)$ and $\varphi : X \rightarrow f_{\ast} Y$ a morphism in $\mathcal{V}(\Shv(D))$. If $(f, \varphi) : (D, X) \rightarrow (E, Y )$ and $(g, \gamma) : (C, W ) \rightarrow (D, X)$ are maps in $\Rep$ then we define the \emph{composite} $(g, \gamma)(f, \varphi): (C, W ) \rightarrow (E,Y)$ as the pair $(fg, (g_{\ast}\varphi)\gamma)$. From the functoriality of $f_{∗}$ and the fact of $(fg)_{\ast} = g_{\ast} f_{\ast}$, it follows that composition in $\Rep$ is well defined and is associative. Moreover, for every $D$ in $\dLat(\Set)$, the identity morphism $id_{D}$ (as a morphism of sites) induces the identity morphism in $\Shv(D)$ so it easily follows that for every pair $(D,X)$ in $\Rep$, $id_{(D,X)}=(id_{D},id_{X})$.
\\

For each morphism $(D,X)$ in $\Rep$ we define $\Gamma(D,X)$ as $X(1)$, and for every $(f,\varphi):(D, X) \rightarrow (E, Y )$ in $\Rep$, define $\Gamma(f,\varphi)=\varphi_{1}:X(1)\rightarrow Y(f(1))=Y(1)$. It easily follows that $\Gamma: \Rep \rightarrow \mathcal{V}$ is a functor.

\section{The representation of \V-models}

In this section we prove that every algebra with RexDFC and CSC in $\Set$ can be represented as an object of the category $\Rep$.
\\

Let $\V$ a variety with RexDFC and CSC. If $A$ is in $\V$ and $\vec{e}\in Z(A)$, recall that from Lemmas \ref{Existential implies Positive} and \ref{Positive implies principal congruence} we get that $\theta^{A}_{\vec{1},\vec{e}}=\theta^{A}(\vec{1},\vec{e})$.

\lemma\label{operations with principal congruences}
Let $A$ be in $\V$ and $\vec{e},\vec{f}\in Z(A)$. The following holds:
\begin{enumerate}
\item $\theta^{A}(\vec{1},\vec{e}\wedge_{A} \vec{f})=\theta^{A}(\vec{1},\vec{e})\vee \theta^{A}(\vec{1},\vec{f})$.
\item $\theta^{A}(\vec{1},\vec{e}\vee_{A} \vec{f})=\theta^{A}(\vec{1},\vec{e})\cap \theta^{A}(\vec{1},\vec{f})$.
\end{enumerate}
\endlemma 
\pf
We prove $1.$ By definition (Subsection \ref{Generalities about Varieties with DFC}), it is clear that $\theta^{A}(\vec{1},\vec{e}\wedge_{A} \vec{f}) \subseteq \theta^{A}(\vec{1},\vec{e})\vee \theta^{A}(\vec{1},\vec{f})$. On the other hand, since  $\vec{e}\wedge_{A} \vec{f}\leq_{A} \vec{e}, \vec{f}$ thus $\theta^{A}(\vec{1},\vec{e}), \theta^{A}(\vec{1},\vec{f})\subseteq \theta^{A}(\vec{1},\vec{e}\wedge_{A} \vec{f})$, hence $\theta^{A}(\vec{1},\vec{e})\vee \theta^{A}(\vec{1},\vec{f})\subseteq \theta^{A}(\vec{1},\vec{e}\wedge_{A} \vec{f})$. The proof of $2.$ is similar.
\epf

\rem\label{Antiiso between Center and FactorCongruences}
Since $\theta^{A}(\vec{1},\vec{1})=\Delta^{A}$ and $\theta^{A}(\vec{1},\vec{0})=\nabla^{A}$, as a direct application of Lemma \ref{operations with principal congruences} it follows that   the map $\phi:Z(A)^{op}\rightarrow FC(A)$ defined by $\phi(\vec{e})=\theta^{A}(\vec{1},\vec{e})$ is an iso of Boolean algebras. 
\endrem

\lemma\label{Technical result 1}
Let $\V$ be a variety with DFC and $A$ an algebra of $\V$. If $\theta\diamond \delta$ in $Con(A)$, then for every $\vec{e}\in A$, the following are equivalent:
\begin{enumerate}
\item $\vec{e}\in Z(A)$
\item $\vec{e}/\theta\in Z(A/\theta)$ and $\vec{e}/\delta\in Z(A/\delta)$.
\end{enumerate}
\endlemma
\pf
Let us assume $\theta\diamond \delta$ in $Con(A)$ and suppose $\vec{e}\in Z(A)$. Without loss of generality we can assume $\vec{e}=(\vec{0},\vec{1})$ in $A=A_{1}\times A_{2}$. Since $DFC$ implies $BFC$ (see the Introduction), from Lemma \ref{FC Factors equivalent BFC} there exist $\alpha_{i}\in Con(A_{i})$, $(i=1,2)$ such that $\theta=\alpha_{1}\times \alpha_{2}$. Thereby, via the canonical isomorphism between $A/\theta$ and $A/\alpha_{1}\times A/\alpha_{2}$ we can conclude that $\vec{e}/\theta=(\vec{0}/\alpha_{1},\vec{1}/\alpha_{2})$, so $\vec{e}/\theta \in Z(A/\theta)$. The proof for $\vec{e}/\delta \in Z(A/\delta)$ is analogue. On the other hand, if $\vec{e}\in Z(A/\theta)$ and $\vec{e}\in Z(A/\delta)$, there exist $A_{1},A_{2}, B_{1}, B_{2}\in \V$ and isomorphisms $\tau_{\theta}:A/\theta\rightarrow A_{1}\times A_{2}$, $\tau_{\delta}:A/\theta\rightarrow B_{1}\times B_{2}$, such that $\tau_{\theta}(\vec{e}/\theta)=(\vec{0}_{A_{1}},\vec{1}_{A_{2}})$ and $\tau_{\delta}(\vec{e}/\delta)=(\vec{0}_{B_{1}},\vec{1}_{B_{2}})$. Since $\theta\diamond \delta$ by assumption, then $A\cong A/\theta \times A/\delta$, so, since $(A_{1}\times A_{2})\times (B_{1}\times B_{2})\cong (A_{1}\times B_{1})\times (A_{2}\times B_{2})=C$, if we write $C_{1}=A_{1}\times B_{1}$ and $C_{2}=A_{2}\times B_{2}$ there exists an isomorphism $\kappa: A\rightarrow C_{1}\times C_{2}$ such that $\kappa(\vec{e})=(\vec{0},\vec{1})$. This concludes the proof.
\epf

\lem\label{Technical result 2}
Let $\V$ be a variety with $BFC$ and $A$ and algebra of $\V$. If $\theta \in FC(A)$ and $\vec{z}/\theta \in Z(A/\theta)$, then there exists an $\vec{e}\in Z(A)$ such that $\vec{e}/\theta=\vec{z}/\theta$.
\endlem
\pf
Let $\delta$ be the factor congruence complementary to $\theta$. Since $\nabla^{A}=\theta\circ \delta$ and $(\vec{z},\vec{1})\in \nabla^{A}$, then, there exists an $e\in A$ such that $(\vec{z},\vec{e})\in \theta$ and $(\vec{e},\vec{1})\in \delta$. It is clear that $\vec{z}/\theta=\vec{e}/\theta\in Z(A/\theta)$ and $\vec{e}/\delta=\vec{1}/\delta \in Z(A/\delta)$. Hence, by Lemma \ref{Technical result 1} we conclude that $\vec{e}\in Z(A)$.
\epf

\lem\label{The presheaf associated a Vmodel}
Let $\V$ be a variety with RexDFC and CSC. For every $A$ algebra of $\V$ and every $\vec{e},\vec{f}\in Z(A)$, if $\vec{f}\leq_{A}\vec{e}$ there exists a (necessarily unique map) $A/\theta_{\vec{1},\vec{e}}\rightarrow A/\theta_{\vec{1},\vec{f}}$ such that the diagram below
\begin{displaymath}
\xymatrix{
A \ar[r] \ar[dr] & A/\theta_{\vec{1},\vec{e}} \ar[d]
\\
 & A/\theta_{\vec{1},\vec{f}}
}
\end{displaymath}
commutes, where the horizontal and diagonal arrows are the respective canonical homomorphisms. Thereby, if $\vec{e}=\vec{f}$ then $A/\theta_{\vec{1},\vec{e}}$ is canonically iso to $A/\theta_{\vec{1},\vec{f}}$.
\endlem
\pf
Again, from Lemmas \ref{Existential implies Positive} and \ref{Positive implies principal congruence}, we obtain that for every $\vec{e}\in Z(A)$, $\theta_{\vec{1},\vec{e}}=\theta(\vec{1},\vec{e})$. If $\vec{f}\leq_{A}\vec{e}$, then $\theta(\vec{1},\vec{e})\subseteq\theta(\vec{1},\vec{f})$, so by Lemma \ref{corollary universal property} the result follows.
\epf

As result, the assignment that sends $\vec{e}\in Z(A)$ to $A/\theta_{\vec{1},\vec{e}}$ is well defined so we obtain a functor $Z(A)^{op}\rightarrow \mathcal{V}$. In conclusion, we have obtained a $\V$-model $\overline{A}$ in $\widehat{Z(A)}$. 

\lem\label{The sheaf associated to a Vmodel}
For every $\V$-model in $\Set$ the presheaf $\overline{A}$ in $\widehat{Z(A)}$ is a sheaf (respect to the coherent coverage on the lattice $Z(A)$).
\endlem
\pf
Since $Z(A)$ is a Boolean algebra, to prove the statement it is enough to verify the sheaf condition for binary partitions, but this leads to a reformulation of item 3. in Lemma \ref{basics about systems}.
\epf

\lem\label{the representation associated to a Vmodel}
For every $\V$-model in $\Set$ the pair $(Z(A), \overline{A})$ is an object of $\Rep$.
\endlem
\pf
We use Lemma \ref{connected models in ShvD}. Since $A$ is an algebra with $\vec{0}$ and $\vec{1}$, it follows that for every $\vec{e}\in Z(A)$, $\overline{A}(\vec{e})$ such condition also holds. Observe that, in the case of $\Shv(D)$, the map $\alpha_{\vec{e}}:Z(\downarrow \vec{e})\rightarrow Z(A/\theta(\vec{1},\vec{e}))$, is canonically defined as $\alpha_{\vec{e}}(\vec{f})=f/\theta(\vec{1},\vec{e})$. We verify that $\alpha_{\vec{e}}$ is biyective. If $\alpha_{\vec{e}}(\vec{f})=\alpha_{\vec{e}}(\vec{g})$, then $[\vec{f},\vec{g}]\in \theta(\vec{1},\vec{e})$. From Lemma \ref{operations with principal congruences}, we have $\theta(\vec{1},\vec{e})\vee \theta(\vec{1},\vec{f})=\theta(\vec{1},\vec{e}\wedge_{A}\vec{f})$ so $\theta(\vec{1},\vec{e})\subseteq \theta(\vec{1},\vec{e}\wedge_{A}\vec{f})$. Since $\vec{e}\wedge_{A} \vec{g}=\vec{g}$ and $[\vec{1},\vec{g}]\in \theta(\vec{1},\vec{g})$, from the transitivity of $\theta(\vec{1},\vec{g})$ we obtain that $[\vec{1},\vec{f}]\in \theta(\vec{1},\vec{g})$ so $\vec{g}\leq_{A} \vec{f}$. The verification of $\vec{f}\leq_{A} \vec{g}$ is similar. Thus $\alpha_{e}$ is injective. To check the surjectivity of $\alpha_{\vec{e}}$, let $\vec{f}/\theta(\vec{1},\vec{e})\in Z(A/\theta(\vec{1},\vec{e}))$. From Lemma \ref{Technical result 2}, there exists a $\vec{z}\in Z(A)$, such that $[\vec{f},\vec{z}]\in \theta(\vec{1},\vec{e})$. Since $[\vec{z},\vec{1}]\in \theta(\vec{1},\vec{z})$, then $[\vec{f},\vec{1}]\in \theta(\vec{1},\vec{e})\vee \theta(\vec{1},\vec{z})=\theta(\vec{1},\vec{e}\wedge_{A}\vec{f})$, again by Lemma \ref{operations with principal congruences}. Thus we obtain that $\theta(\vec{f},\vec{1})\subseteq \theta(\vec{1},\vec{e}\wedge_{A}\vec{f})$. From Lemma \ref{Useful lema Centrals}, $[\vec{e}\wedge_{A}\vec{z}, \vec{z}]\in \theta(\vec{1},\vec{e})$, so, since $[\vec{f},\vec{z}]\in \theta(\vec{1},\vec{e})$, we get that $[\vec{f},\vec{e}\wedge_{A}\vec{z}]\in \theta(\vec{1},\vec{e})$. Hence, since $\vec{e}\wedge_{A}\vec{z}\leq_{A} \vec{f}$ and  $\theta(\vec{f},\vec{1})\subseteq \theta(\vec{1},\vec{e}\wedge_{A}\vec{f})$, we conclude that $\alpha_{\vec{e}}(\vec{e}\wedge_{A}\vec{z})=f/\theta(\vec{1},\vec{e})$. Finally, if $A/\theta(\vec{1},\vec{e})$ is trivial, it follows that $[\vec{0},\vec{1}]\in \theta(\vec{1},\vec{e})$, thus $\theta(\vec{1},\vec{0})=\theta(\vec{1},\vec{e})$. Hence, by Remark \ref{Antiiso between Center and FactorCongruences}, $\vec{e}$ must be $\vec{0}$. This concludes the proof.
\epf

\section{RexDFC and CSC induce homomorphisms of Boolean algebras}

As we saw in Section \ref{Coextensivity and Center Stabitlity}, not every variety with BFC has center stable. In this section we prove that a variety with RexDFC having center stable by complements is in fact a variety with the Fraser Horn Property. This result will allow us to prove that the every homomorphism $f$ in the variety induces a Boolean algebra homomorphism between the centers of $dom(f)$ and $cod(f)$. 

\lem[Theorem 1 \cite{FH1970}]\label{Theorem 1 FHP}
Let $\textbf{K}$ be a variety and $A$, $B$ be algebras of $\textbf{K}$. The following are equivalent:
\begin{enumerate}
\item $\textbf{K}$ has FHP.
\item For every $A,B\in \textbf{K}$ and $\gamma\in Con(A\times B)$,
\[\Pi_{1}\cap (\Pi_{2} \vee \gamma) \subseteq \gamma\; \textrm{and}\; \Pi_{2}\cap (\Pi_{1} \vee \gamma) \subseteq \gamma\]
where $\Pi_{1}$ is the kernel of the projection on $A$ and $\Pi_{2}$ is the kernel of the projection on $B$.
\end{enumerate} 
\endlem

\lem[Theorem 3 \cite{FH1970}]\label{Theorem 3 FHP}
Let $A$ and $B$ be similar algebras. The following are equivalent:
\begin{enumerate}
\item $A\times B$ has FHP.
\item For every $a,c\in A$ and $b,d\in B$,
\[\theta^{A\times B}((a,b),(c,d))=\theta^{A}(a,c)\times \theta^{B}(b,d) \] 
\end{enumerate}
\endlem

\lem\label{useful lema Malsev}
Let $A$ and $B$ be algebras with finite $n$-ary function symbols and $f:A\rightarrow B$ an homomophism. If $(a,b)\in \theta^{A}(\vec{c},\vec{d})$, then $(f(a),f(b))\in \theta^{B}(f(\vec{c}),f(\vec{d}))$. Thus, if $[\vec{a},\vec{b}]\in \theta^{A}(\vec{c},\vec{d})$ then $[f(\vec{a}),f(\vec{b})]\in \theta^{A}(f(\vec{c}),f(\vec{d}))$.
\endlem
\pf
Apply Lemma \ref{Gratzer Malsev Lemma}. 
\epf

Let $\V$ be a variety with DCF. As we have seen, for every algebra $A\in \V$ and $\vec{e}\in Z(A)$ the span $A/\theta_{\vec{0},\vec{e}}\leftarrow A \rightarrow A/\theta_{\vec{1},\vec{e}}$ is a product. Notice that in this case $\Pi_{1}=\theta_{\vec{0},\vec{e}}$ and $\Pi_{2}=\theta_{\vec{1},\vec{e}}$.

\lem\label{RexDFC implies FHP}
Let $\V$ be a variety with RexDFC, $A\in \V$ and $\vec{e},\vec{f}\in Z(A)$ such that $\vec{e}\diamond_{A}\vec{f}$. If for every $\gamma\in Con(A)$, $\vec{e}/\gamma \diamond_{A/\gamma} \vec{f}/\gamma$, then $\V$ has the FHP. 
\endlem
\pf
Since $\V$ has RexDFC, from Lemmas \ref{Existential implies Positive} and \ref{Positive implies principal congruence}, then,  for every $\vec{e}\in Z(A)$, $\theta^{A}_{\vec{1},\vec{e}}=\theta^{A}(\vec{1},\vec{e})$. So, if $\vec{e}\diamond_{A}\vec{f}$ then $\theta^{A}_{\vec{0},\vec{e}}=\theta^{A}(\vec{1},\vec{f})$. We use Lemma \ref{Theorem 1 FHP}. To do so, we prove $\theta^{A}(\vec{1},\vec{e})\cap (\theta^{A}(\vec{1},\vec{f})\vee \gamma)\subseteq \gamma$. Suppose $(x,y)\in \theta^{A}(\vec{1},\vec{e})\cap (\theta^{A}(\vec{1},\vec{f})\vee \gamma)\subseteq \gamma$, then, $(x,y)\in \theta^{A}(\vec{1},\vec{e})$ and there are $c_{0},...,c_{N}\in A$, with $c_{0}=x$ and $c_{N}=y$, such that $(c_{2i},c_{2i+1})\in \theta^{A}(\vec{1},\vec{f})$ and $(c_{2i+1},c_{2(i+1)})\in \gamma$. Since $A\rightarrow A/\gamma$ is clearly an homomorphism, from Lemma \ref{useful lema Malsev}, we obtain that $(x/\gamma,y/\gamma)\in \theta^{A/\gamma}(\vec{1}/\gamma,\vec{e}/\gamma)$, $(c_{2i}/\gamma,c_{2i+1}/\gamma)\in \theta^{A/\gamma}(\vec{1}/\gamma,\vec{f}/\gamma)$ and $c_{2i+1}/\gamma=c_{2(i+1)}/\gamma$. From transitivity of $\theta^{A/\gamma}(\vec{1}/\gamma,\vec{f}/\gamma)$, we get that $(x/\gamma,y/\gamma)\in \theta^{A/\gamma}(\vec{1}/\gamma,\vec{f}/\gamma)$. Therefore, $(x/\gamma,y/\gamma)\in \theta^{A/\gamma}(\vec{1}/\gamma,\vec{e}/\gamma) \cap \theta^{A/\gamma}(\vec{1}/\gamma,\vec{f}/\gamma)=\Delta^{A/\gamma}$, since $\vec{e}/\gamma \diamond_{A/\gamma} \vec{f}/\gamma$ by assumption, so $(x,y)\in \gamma$. The proof of $\theta^{A}(\vec{1},\vec{f})\cap (\theta^{A}(\vec{1},\vec{e})\vee \gamma)\subseteq \gamma$ is similar. This concludes the proof.
\epf

\corollary\label{RexDFC with SCC has FHP}
Let $\V$ be a variety with RexDFC. If $\V$ has SCC then has FHP. 
\endcorollary
\pf
Inmmediate from Lemma \ref{RexDFC implies FHP}.
\epf

\lem\label{FHP implies TexDFC}
Every variety $\V$ with FHP is TexDFC.
\endlem
\pf
We want to prove there exists an existential formula $\varphi$ which defines $\theta_{\vec{0},\vec{e}}$ in terms of $\vec{e}$. To do so, let $C\in \V$ and $\vec{e}\in Z(C)$. Let us consider $\varphi(x,y,\vec{z})=\pi(x,y,\vec{0},\vec{z})$, where $\pi(x,y,\vec{0},\vec{z})$ is the formula of Lemma \ref{Malsev restated}. It is clear that $\varphi$ is existential. Let $A,B\in \V$ and $(a,b),(c,d)\in A\times B$. If $A\times B\models \varphi((a,b),(c,d), (\vec{0},\vec{1}))$, then from Lemma \ref{Malsev restated}, $((a,b),(c,d))\in \theta^{A\times B}((\vec{0},\vec{0}),(\vec{0},\vec{1}))$. Since $\V$ has FHP by hypothesis, then from Lemma \ref{Theorem 3 FHP} $\theta^{A\times B}((\vec{0},\vec{0}),(\vec{0},\vec{1}))=\theta^{A}(\vec{0},\vec{0})\times \theta^{B}(\vec{0},\vec{1})=\Delta^{A}\times \nabla^{B}$. Hence, $a=c$. On the other hand, suppose $a=c$. Let $P=\textbf{F}(x)\times \textbf{F}(x,y)$, and consider the pair $((x,x),(x,y))\in P$. Since $\V$ has FHP by assumption, thus, again by Lemma \ref{Theorem 3 FHP}, $\theta^{P}((\vec{0},\vec{0}),(\vec{0},\vec{1}))=\theta^{\textbf{F}(x)}(\vec{0},\vec{0})\times \theta^{\textbf{F}(x,y)}(\vec{0},\vec{1})=\Delta^{\textbf{F}(x)}\times \nabla^{\textbf{F}(x,y)}$. Observe that $((x,x),(x,y))\in \theta^{P}((\vec{0},\vec{0}),(\vec{0},\vec{1}))$ and consider the assignments $\alpha_{A}:\{x\}\rightarrow A$ and $\alpha_{B}:\{x,y\}\rightarrow B$, defined as $\alpha_{A}(x)=a$ and $\alpha_{B}(x)=b$, $\alpha_{B}(y)=d$, respectively. From the left adjointness of the free functor $\textbf{F}: \Set \rightarrow \mathcal{V}$ to the forgetful functor, there exist a unique pair of homomorphisms $\beta_{A}:\textbf{F}(x)\rightarrow A$ and $\beta_{B}:\textbf{F}(x,y)\rightarrow B$ extending $\alpha_{A}$ and $\alpha_{B}$. Consider the morphism $g=\beta_{A}\times \beta_{B}: P\rightarrow A\times B$. From Lemma \ref{useful lema Malsev} we obtain that $(g(x,x),g(x,y))=((a,b),(a,c))\in \theta^{A\times B}((\vec{0},\vec{0}),(\vec{0},\vec{1}))$. Hence, by Lemma \ref{Malsev restated}, $A\times B\models \varphi((a,b),(c,d), (\vec{0},\vec{1}))$. Therefore, $\V$ has LexDFC. The proof for $\V$ has RexDFC is analogue. This concludes the proof.
\epf

As a straight consequence of Corollary \ref{RexDFC with SCC has FHP} and Lemma \ref{FHP implies TexDFC} we obtain

\corollary\label{RexDFC and SCC implies TexDFC}
Every variety $\V$ with RexDFC and CSC is TexDFC.
\endcorollary

\lem\label{RexDFC and CS entails Lattice morphisms}
Let $\V$ be a variety with RexDFC, $A,B\in \V$ and $f:A\rightarrow B$ be an homomorphism. If $f$ preserves central elements, then $f|_{Z(A)}:Z(A)\rightarrow Z(B)$ is a bounded lattice homomorphism.
\endlem
\pf
First of all, observe that from Corollary \ref{RexDFC and SCC implies TexDFC}; and Lemmas \ref{Existential implies Positive} and \ref{Positive implies principal congruence}, we get that for every $\vec{e}\in A$, $\theta^{A}_{\vec{0},\vec{e}}=\theta^{A}(\vec{0},\vec{e})$ and $\theta^{A}_{\vec{1},\vec{e}}=\theta^{A}(\vec{1},\vec{e})$. Now, since $f$ is homomorphism, it is clear that preserves $\vec{0}$ and $\vec{1}$. So, if $\vec{e}_{1},\vec{e}_{2}\in Z(A)$, and $\vec{a}=\vec{e}_{1}\wedge_{A}\vec{e}_{2}$, thus from Lemma \ref{Useful lema Centrals}, $[\vec{0},\vec{a}]\in \theta^{A}(\vec{0},\vec{e}_{1})$ and $[\vec{a},\vec{e}_{2}]\in \theta^{A}(\vec{1},\vec{e}_{1})$. Thus, since $f(\vec{e})\in Z(A)$ for every $\vec{e}\in A$ by hypothesis; from Lemma \ref{useful lema Malsev} we get that $[\vec{0},f(\vec{a})]\in \theta^{A}(\vec{0},f(\vec{e}_{1}))$ and $[f(\vec{a}),f(\vec{e}_{2})]\in \theta^{B}(\vec{1},f(\vec{e}_{1}))$ so again by Lemma \ref{Useful lema Centrals} we can conclude that $f(a)=f(\vec{e}_{1})\wedge_{B}f(\vec{e}_{2})$. The proof for the preservation of the join is similar. 
\epf

A direct application of Lemma \ref{RexDFC and CS entails Lattice morphisms} gives as result

\corollary\label{RexDFC with CSC entails Boolean morhisms} 
Let $\V$ be a variety with RexDFC and SCC. Then, for every $A,B\in \V$ and every homomorphism $f:A\rightarrow B$, the map $f|_{Z(A)}:Z(A)\rightarrow Z(B)$ is an homomorphism of Boolean algebras. 
\endcorollary

\section{The representation theorem}

For the rest of this section $\V$ will be a variety with RexDFC and CSC. Next we show that the functor $\Gamma : \Rep \rightarrow \mathcal{V}$ has a fully faithful left adjoint.
\\

Let $A$ and $B$ be $\V$-models in $\Set$ and let $f: A\longrightarrow B$ be a ${\mathcal{V}}$. Since $\V$ is CSC, from Corollary \ref{RexDFC with CSC entails Boolean morhisms}, the restriction of $f$ to $Z(A)$ determines a morphism of boolean algebras $f:Z(A)\rightarrow Z(B)$. Such morphism is also a morphism of lattices so determines a geometric morphism $f:\Shv(Z(B)))\rightarrow \Shv(Z(B))$ whose direct image $f_{\ast}$ is defined as $f_{\ast}(G)(\vec{e})=G(f(\vec{e}))$, for every $G\in \Shv(Z(B))$.

\lem\label{V morphisms determines natural transformations}
Every morphism $f:A\rightarrow B$ in $\mathcal{V}$ determines a natural transformation $\overline{f}: \overline{A} \longrightarrow f_{\ast}(\overline{B})$ in $\Shv(Z(A))$.
\endlem

\pf
Let $\vec{e}\in Z(A)$. If $i_{\vec{e}}:A\rightarrow A/\theta^{A}(\vec{1},\vec{e})$ and $i_{f(\vec{e})}:B\rightarrow A/\theta^{A}(\vec{1},f(\vec{e}))$ are the canonical homomorphisms, from Corollary \ref{corollary universal property}

\begin{displaymath}
\xymatrix{
A\ar[r]^-{f} \ar[d]_-{i_{\vec{e}}} & B\ar[d]^-{i_{f(\vec{e})}}
\\
A/\theta^{A}(\vec{1},\vec{e})\ar[r]_-{f_{e}} & B/\theta^{B}(\vec{1},f(\vec{e}))
}
\end{displaymath}

it follows that there exists a unique morphism $f_{\vec{e}}: A/\theta^{A}(\vec{1},\vec{e})\rightarrow B/\theta^{B}(\vec{1},f(\vec{e}))$ in $\mathcal{V}$, such that the diagram above commutes. Consider the assigment $\overline{f}: \overline{A} \longrightarrow f_{\ast}(\overline{B})$, defined as $\overline{f}_{\vec{e}}=f_{\vec{e}}$. We prove that $\overline{f}$ is natural in $\Shv(Z(A))$. Let $\vec{e}_1,\vec{e}_2\in Z(A)$ with $\vec{e}_2\leq_{A} \vec{e}_1$. From Lemma \ref{RexDFC and CS entails Lattice morphisms} it follows that $f(\vec{e}_2)\leq_{B} f(\vec{e}_1)$, so again by Corollary \ref{corollary universal property}, the diagram below (where the rows of the right square are the canonical morphisms $A/\theta^{A}(\vec{1},\vec{e}_{1})\rightarrow A/\theta^{A}(\vec{1},\vec{e}_{2})$ and $B/\theta^{B}(\vec{1},f(\vec{e}_{1}))\rightarrow B/\theta^{B}(\vec{1},f(\vec{e}_{2}))$, respectively), 

\begin{displaymath}
\xymatrix{
A \ar[d]_-{f} \ar[r]^-{i_{\vec{e}_{1}}} & A/\theta^{A}(\vec{1},\vec{e}_{2}) \ar[r] \ar[d]^-{f_{\vec{e}_{1}}} & A/\theta^{A}(\vec{1},\vec{e}_{1}) \ar[d]^-{f_{\vec{e}_{2}}}
\\
B \ar[r]_-{i_{f(\vec{e}_{1})}} & B/\theta^{B}(\vec{1},f(\vec{e}_{1})) \ar[r] & B/\theta^{B}(\vec{1},f(\vec{e}_{2}))
}
\end{displaymath}

commutes. Since $ B/\theta^{B}(\vec{1},f(\vec{e}_{1}))=f_{\ast}(\overline{B})(\vec{e})$, the result follows.
\epf

Lemmas \ref{the representation associated to a Vmodel} and \ref{V morphisms determines natural transformations} allows us to define an assigment $\mathcal{F}:\mathcal{V}\rightarrow \Rep$ as $\mathcal{F}(A)=(Z(A), \overline{A})$ and $\mathcal{F}(f:A\rightarrow B)=(f,\overline{f})$.

\lem\label{The assignment is functorial}
The assignment $\mathcal{F}:\mathcal{V}\rightarrow \Rep$ is functorial. 
\endlem
\pf
It is clear that $\mathcal{F}(id_{A})=(id_{A}, id_{\overline{A}})=id_{\mathcal{F}(A)}$. So, let $f:A\rightarrow B$ and $h:B\rightarrow C$ be morphisms in ${\mathcal{V}}$. Then, we get that $\mathcal{F}(hf)=(hf,\overline{hf})$ and $\mathcal{F}(h)\mathcal{F}(f)=(hf,f_{\ast}(\overline{h})\overline{f}))$. 
\begin{displaymath}
\xymatrix{
A \ar[r]^-{f} \ar[d]_-{i_{\vec{e}}} & B \ar[r]^-{h} \ar[d]_-{i_{f(\vec{e})}} & C \ar[d]^-{i_{hf(\vec{e})}}
\\
A/\theta^{A}(\vec{1},\vec{e}) \ar[r]_-{f_{\vec{e}}} \ar@/_2pc/[rr]_-{hf_{\vec{e}}} & B/\theta^{B}(\vec{1},f(\vec{e})) \ar[r]_-{h_{f(\vec{e})}} & C/\theta^{C}(\vec{1},hf(\vec{e}))
}
\end{displaymath}

Since  $i_{hf(\vec{e})}=i_{h(f(\vec{e}))}$ and ${hf(\vec{e})}={h(f(\vec{e}))}$, then from  Corollary \ref{corollary universal property}, the diagram above commutes for every $\vec{e}\in Z(A)$. Hence, since $\overline{h}_{f(\vec{e})}\overline{f}_{\vec{e}}=f_{\ast}(\overline{h})_{\vec{e}}\overline{f}_{\vec{e}}$, then $f_{\ast}(\overline{h})\overline{f}=\overline{hf}$. Thereby $\mathcal{F}(h)\mathcal{F}(f)=\mathcal{F}(hf)$.
\epf

\lem\label{Previous TheAdjunction 1}
Let $A$ be an algebra of $\V$ in $\Set$ and $P$ be a $\V$-model in $\widehat{Z(A)}$. For every homomorphism $g:A\rightarrow P(\vec{1})$ in $\V$, the following are equivalent:
\begin{enumerate}
\item For every $\vec{e}\in Z(A)$, $g(\vec{e})\cdot \vec{e}=\vec{1} \in P(\vec{1})$,
\item There exist a unique morphism of $\V$-models $\phi: \overline{A}\rightarrow P$ in $\widehat{Z(A)}$, such that $\phi_{\vec{1}}=g$.
\end{enumerate}
\endlem
\pf
Let us assume $g(\vec{e})\cdot \vec{e}=\vec{1} \in P(\vec{1})$, for every $\vec{e}\in Z(A)$. Since $\vec{1}\cdot \vec{e}=\vec{1}$ for every $\vec{e}\in Z(A)$, from the universal property of $A\rightarrow A/\theta^{A}(\vec{1},\vec{e})$ (Corollary \ref{corollary universal property}), for every $\vec{e}\in Z(A)$, there exists a unique homomorphism $A/\theta^{A}(\vec{1},\vec{e})\rightarrow P(\vec{e})$ in $\V$ such that the diagram below
\begin{displaymath}
\xymatrix{
A\cong \overline{A}(\vec{1}) \ar[d] \ar[r] \ar@/^1pc/[r]^-{g} & P(1) \ar[d]
\\
A/\theta^{A}(\vec{1},\vec{e})=\overline{A}(\vec{e}) \ar[r]_-{\phi_{\vec{e}}} & P(\vec{e})
}
\end{displaymath}

commutes. Observe that Corollary \ref{corollary universal property} also grants that the collection $\{\phi_{\vec{e}}\mid \vec{e}\in Z(A)\}$ is natural. The proof of the last part follows from the naturality of $\phi$.
\epf

\lem\label{Previous TheAdjunction 2}
Let $A$ be an algebra of $\V$ in $\Set$ and $(E,Y)$ in $\Rep$. For every $g:A\rightarrow Y(1)$ in $\V$, the following are equivalent:
\begin{enumerate}
\item There is a unique lattice morphism $f:Z(A)\rightarrow E$, such that, for every $\vec{e}\in Z(A)$, $g(\vec{e})\cdot f(\vec{e})=\vec{1}\in f_{\ast}(Y)(\vec{e})$.
\item There exists a unique $(f,\varphi):(Z(A),\overline{A})\rightarrow (E,Y)$ in $\Rep$, such that $\phi_{\vec{1}}=g$.
\end{enumerate}
\endlem
\pf
If we assume $2.$ then $1.$ is granted for the naturality of $\phi$. On the other hand, by assuming $1.$, it follows that $f_{\ast}(Y)$ is in $\Shv(Z(A))$, so from Lemma \ref{Previous TheAdjunction 1}, for the map $g:A\rightarrow Y(f(\vec{1}))=f_{\ast}(Y)(\vec{1})$ there exists a unique morphism of $\V$-models $\phi:\overline{A}\rightarrow f_{\ast}(Y)$ such that $\phi_{\vec{1}}=g$. Thereby, the uniqueness of $(f,\varphi):(Z(A), \overline{A})\rightarrow (E,Y)$ in $\Rep$, easily follows. This concludes the proof.
\epf

Coming up next, we prove the main result of this paper.

\theorem\label{TheAdjunction}
The functor $\Gamma:\Rep \rightarrow \mathcal{V}$ has a full and faithful left adjoint.
\endtheorem
\pf
Let $A$ be an arbitrary algebra of $\V$. From Lemma \ref{the representation associated to a Vmodel}, $(Z(A),\overline{A})$ is an object of $\Rep$. Let us to consider the (iso) map $A\rightarrow A/\theta^{A}(\vec{1},\vec{1})=\overline{A}(\vec{1})=\Gamma(Z(A),\overline{A})$. We prove this map is universal from $A$ to $\Gamma$. To do so, let $(C,X)$ be in $\Rep$ and $g:A\rightarrow X(1)=\Gamma(C,X)$ be an arbitrary morphism of $\mathcal{V}$. From the center stability of $\V$, for every $\vec{e}\in Z(A)$, $g(\vec{e})\in Z(X(1))$. Since $(C,X)$ is in $\Rep$, $X$ is connected in $\Shv(C)$, so by Proposition \ref{connected models in ShvD}, there are bijections $\alpha_{1},\alpha_{g(\vec{e})}$ making the diagram below
\begin{displaymath}
\xymatrix{
Z(X(1)) \ar[r]^-{\alpha_{1}} \ar[d]_-{k_{g(\vec{e})}} & Z(C) \ar[d]^-{j_{g(\vec{e})}}
\\
Z(X(g(\vec{e}))) \ar[r]_-{\alpha_{g(\vec{e})}} & Z(\downarrow g(\vec{e}))
}
\end{displaymath}

commutes, for every $\vec{e}\in Z(A)$ (with $k_{g(\vec{e})}(\vec{h})=\vec{h}\cdot g(\vec{e})$ and $j_{g(\vec{e})}(l)=l\wedge g(\vec{e})$). Let us define $f:Z(A)\rightarrow C$, as $f(\vec{e})=\alpha_{\vec{1}}(g(\vec{e}))$. By Lemma \ref{RexDFC and CS entails Lattice morphisms}, $f$ is a  lattice morphism, thus $f_{\ast}(Y)$ is in $\Shv(Z(A))$. The commutativity of diagram above allows us to make the following calculation

\[\alpha_{f(\vec{e})}(g(\vec{e})\cdot f(\vec{e}))=\alpha_{\vec{1}}(g(\vec{e}))\wedge f(\vec{e})= f(\vec{e}) = \alpha_{f(\vec{e})}(\vec{1}\cdot f(\vec{e}))\]

Hence, $g(\vec{e})\cdot f(\vec{e})=\vec{1} \in f_{\ast}(Y)(\vec{e})$. Thereby, from Lemma \ref{Previous TheAdjunction 1}, there exists a unique morphism of $\V$-models $\phi: \overline{A}\rightarrow f_{\ast}(Y)$ in $\Shv(Z(A))$, such that $\phi_{\vec{1}}=g$, so, by Lemma \ref{Previous TheAdjunction 2} there exist a unique $(f,\varphi):(Z(A),\overline{A})\rightarrow (E,Y)$ in $\Rep$, such that $\phi_{\vec{1}}=g$.

\section{Corollaries in terms of local homeos}\label{Corollaries in terms of local homeos}

It is a classical result that for any topological space $X$, the category ${\LH{X}}$ of local homeomorphisms over $X$ is equivalent to the topos $\Shv(X)$ of sheaves over the same space (see Section II.6 in \cite{MM2012}). The equivalence $\Shv(X) \rightarrow \LH{X}$ sends a sheaf $P : \mathcal{O}(X) \rightarrow \Set$ to the \emph{bundle of germs} of $P$ defined as follows. For each $x \in X$, let $P_{x} = \limit_{x\in U} P(U)$ where the colimit is taken over the poset of open neighborhoods of $x$ (ordered by reverse inclusion). The family of $P_{x}$'s determines a function $\pi :\sum_{x\in X} P_{x} \rightarrow X$. Also, each $s \in P(U)$ determines an obvious function $\dot{s} : U \rightarrow \sum_{x\in X} P_{x}$ such that
$\pi \dot{s} : U \rightarrow X$ is the inclusion $U \rightarrow X$. The set $\sum_{x\in X} P_{x}$ is topologized by taking as a base of opens all the images of the functions $\dot{s}$. This topology makes $\pi$ into a local homeo, the above mentioned bundle of germs.
\\

Any basis B for the topology of X may be considered as a subposet $B \rightarrow \mathcal{O}(X)$. The usual Grothendieck topology on $\mathcal{O}(X)$ restricts along $B \rightarrow \mathcal{O}(X)$ and the resulting morphism of sites determines an equivalence $\Shv(B) \rightarrow \Shv(X)$; see Theorem II.1.3 in \cite{MM2012}. The composite equivalence $\Shv(B) \rightarrow \Shv(X) \rightarrow \LH{X}$ is very similar to the previous one because, by finality (in the sense of Section IX.3 of \cite{M1971}), the colimit $P_{x} = \limit_{x\in U} P(U)$ may be calculated using only basic open sets.

According to \cite{S1980}, the \emph{spectrum} of a distributive lattice $D$ is the topological space $\sigma D$ whose points are the lattice morphisms $D \rightarrow \textbf{2}$ and whose topology has, as a basis, the subsets $\sigma(a) \subseteq \sigma D$ (with $a \in D$) defined by $\sigma(a) = \{p \in \sigma D \mid p(a) = 1 \in \textbf{2}\} \subseteq \sigma D$. In this way, we may identify $D$ with the basis of its spectrum and obtain an equivalence $\Shv(D) \rightarrow \LH{\sigma D}$. It assigns to each sheaf $P : D^{op} \rightarrow \Set$ the local homeomorphism whose fiber $P_{p}$ over the point $p : D \rightarrow \textbf{2}$ in $\sigma D$ is $P_{p} = \limit_{p\in \sigma(a)} P(a)$. 
\\

Let $\textbf{V}$ be a variety with RexDFC, $A$ be an algebra of $\textbf{V}$ and consider its center $Z(A)$. From the formulation of above, it can be proved that the points of $\sigma(Z(A))$ can be identified with the ultrafilters of $Z(A)$ and the basis $\{\sigma (\vec{e})\mid \vec{e}\in Z(A)\}$ becomes a basis of clopens, making the space $\sigma (Z(A))$ a \emph{Stone space} (\cite{J1982}). 

This facts, together with the ones considered before, allows us to obtain an equivalence $\Shv(Z(A)) \rightarrow \LH{\sigma Z(A)}$ that sends a sheaf $P \in \Shv(Z(A))$ to a local homeomorphism over $\sigma Z(A)$, whose fiber $P_{U}$ over an ultrafilter $U$ in $\sigma Z(A)$ may be described as \[P_{U}=\limit_{\vec{e}\in U}P(\vec{e})\]

\lem\label{Fibers of the sheaf}
Let $\textbf{V}$ be a variety with RexDFC, $A$ be an algebra of $\textbf{V}$. For every ultrafilter $U$ of $Z(A)$, there exists a unique isomorphism $A/\theta(U)\rightarrow \limit_{\vec{e}\in U^{op}} A/\theta(\vec{1},\vec{e})$ such that the following diagram
\begin{displaymath}
\xymatrix{
A \ar[r] \ar[d] & A/\theta(U) \ar[d] 
\\
A/\theta(\vec{1},\vec{e}) \ar[r] & \limit_{\vec{e}\in U^{op}} A/\theta(\vec{1},\vec{e})
}
\end{displaymath}
commutes for every $\vec{e}\in Z(A)$.
\endlem
\pf
Let $\vec{e}\in Z(A)$ and $U$ be an ultrafilter of $Z(A)$. If $\vec{e}\in U$, then $[\vec{1},\vec{e}]\in \theta(U)$, so $\theta(\vec{1},\vec{e})\subseteq \theta(F)$. Thus, from Corollary \ref{corollary universal property}, there exists a unique homomorphism $\rho_{\vec{e}}: A/\theta(\vec{1},\vec{e})\rightarrow A/\theta(U)$ such that the diagram 
\begin{displaymath}
\xymatrix{
A \ar[r]^-{\nu_{\vec{e}}} \ar[dr]_-{\nu_{F}} & A/\theta(\vec{1},\vec{e}) \ar[d]^-{\rho_{\vec{e}}} 
\\
 & A/\theta(U) 
}
\end{displaymath}

We prove that the sink $H=\{\rho_{\vec{e}}\mid \vec{e}\in U\}$ is a colimit. Let $\vec{e},\vec{f}\in U$, such that $\vec{e}\leq_{A}\vec{f}$, then $\vec{e}\wedge_{A}\vec{f}=\vec{e}$. From Lemma \ref{Useful lema Centrals}, $[\vec{e},\vec{f}]\in \theta(\vec{1},\vec{e})$, then, since $[\vec{1},\vec{e}]\in \theta(\vec{1},\vec{e})$, we get that $[\vec{1},\vec{f}]\in \theta(\vec{1},\vec{e})$ and consequently $\theta(\vec{1},\vec{f})\subseteq \theta(\vec{1},\vec{e})$. This proves that $H$ is natural. In order to verfify that $H$ is a cocone, let $\nu_{\vec{e}}:A\rightarrow A/\theta(\vec{1},\vec{e})$, $\nu_{\vec{f}}:A\rightarrow A/\theta(\vec{1},\vec{f})$ and $\nu_{U}:A\rightarrow A/\theta(U)$ be the canonical homomorphisms. From the diagram below
\begin{displaymath}
\xymatrix{
A \ar[r]_-{\nu_{\vec{f}}} \ar[dr]_-{\nu_{\vec{e}}} \ar@/^1pc/[rr]^{\nu_{U}} & A/\theta(\vec{1},\vec{f}) \ar[d]^-{\lambda} \ar[r]_-{\rho_{\vec{e}}} & A/\theta(U)  
\\
 & A/\theta(\vec{1},\vec{e}) \ar[ur]_-{\rho_{\vec{e}}} & 
}
\end{displaymath}

we obtain that $\nu_{F}=\rho_{\vec{f}}\nu_{\vec{f}}=\rho_{\vec{e}}\nu_{\vec{e}}$ and $\lambda \nu_{\vec{f}}=\nu_{\vec{e}}$. Then, from the following calculations 
\[\rho_{\vec{f}}\nu_{\vec{f}}= \rho_{\vec{e}}\nu_{\vec{e}}=\rho_{\vec{e}}\lambda \nu_{\vec{f}}\]
and the fact of $\nu_{\vec{f}}$ is epi, we conclude $\rho_{\vec{f}}=\rho_{\vec{e}}\lambda$. Finally, we check that $H$ is universal. Let us to consider the commutative diagram

\begin{displaymath}
\xymatrix{
A/\theta(\vec{1},\vec{f}) \ar[r]^-{\alpha_{\vec{f}}} \ar[d]_-{\lambda} & B 
\\
A/\theta(\vec{1},\vec{e}) \ar[ur]_-{\alpha_{\vec{e}}}  & 
}
\end{displaymath}

From the commutativity of the left diagram below, 

\begin{center}
\begin{tabular}{ccc}

\xymatrix{ 
A \ar@/^1.5pc/[rr]^{\nu_{U}} \ar[r]_-{\nu_{\vec{f}}} \ar[dr]_-{\nu_{\vec{e}}} & A/\theta(\vec{1},\vec{f}) \ar[r]^-{\rho_{\vec{f}}} \ar[d]_-{\lambda} \ar[dr]^-{\alpha_{\vec{f}}} & A/\theta(U) \ar[d]^-{\mu}
\\
& A/\theta(\vec{1},\vec{e}) \ar[r]_-{\alpha_{\vec{e}}}  & B
}

& & \xymatrix{
A \ar[r]^-{\nu_{U}} \ar[dr]_-{g} & A/\theta(U) \ar[d]^-{\mu} 
\\
 & B 
}
\end{tabular}
\end{center}

we can deduce that the morphism $g=\alpha_{\vec{e}}\nu_{\vec{e}}$ identifies the elements of $U$, then, by Lemma \ref{universal property principal congruences}, there exists a unique morphism $\mu: A/\theta(U)\rightarrow B$, such that the right upper diagram cummutes. From Corollary \ref{corollary universal property}, for every $\vec{e}\in Z(A)$ we can conclude that the diagram  
\begin{displaymath}
\xymatrix{
A/\theta(\vec{1},\vec{e}) \ar[r]^-{\rho_{\vec{e}}} \ar[dr]_-{\alpha_{\vec{e}}} & A/\theta(U) \ar[d]^-{\mu} 
\\
 & B 
}
\end{displaymath}
The result follows from the universal property of colimits.
\epf

Thus, if $\V$ is a variety with RexDFC and $A\in \V$, from Lemma \ref{Fibers of the sheaf} and the discussion along this section, in the particular case of the representing sheaf $\overline{A}$ in $\Shv(Z(A))$, the fiber over an ultrafilter $U$ of $Z(A)$ is 
\[\overline{A}_{U}=\limit_{\vec{e}\in U}\overline{A}(\vec{e})=\limit_{\vec{e}\in U} A/\theta(\vec{1},\vec{e})= A/\theta(U)\] 
That is, if we consider the representing sheaf  $\overline{A}$ in $\Shv(Z(A))$ as a local homemorphism over $\sigma Z(A)$, then the fiber over a point $U$ in $\sigma Z(A)$ is the quotient of $A$ by the principal congruence of $A$ containing $U$. 
\\

Recall that a topos $\mathcal{E}$ with subobject classifier $\Omega$ is boolean, if the cospan $\top,\perp: 1\rightarrow \Omega$ is a coproduct (or equivalently $\Omega \cong 1+1$). Thereby, if $\V$ is a variety with RexDFC, and $A\in V$, from Lemma \ref{LemCoprosInShvD}, it turns that $\Shv(Z(A))$ is boolean a topos, and, in particular, a boolean coherent category. From Lemma \ref{Definability of beeing factor} and Definition \ref{Connected VModels in a topos}, it follows that the theory of internal connected $\Vmod$-models $\mathbb{C}_{\Vmod}$ is a first order theory, so if we call $\mathbb{C}_{\Vmod}'$ to the Morleyzation of $\mathbb{C}_{\Vmod}$, from D1.5.13 of \cite{J2002}, we obtain that $\mathbb{C}_{\Vmod}(\Shv(Z(A)))\simeq \mathbb{C}_{\Vmod}'(\Shv(Z(A)))$.

It is known that every point $x:1\rightarrow X$ of a topological space $X$ determines a geometric morphism $\Set \rightarrow \LH{X}$ whose inverse image $\LH{X}\rightarrow \Set$ sends a local homeomorphism to the corresponding fiber over $x$. Since geometric morphisms preserve the interpretation of coherent sequents, they preserve internal connected $\Vmod$-models.	As as a consequence of the latter discussion, we obtain the following result.

\corollary\label{The fibers of Pierce Sheaf are connected}
Let $\V$ be a variety with RexDFC and CSC. Then, every algebra of $\V$ can be represented as the algebra of global sections of a local homeomorphism (over the Stone space $\sigma Z(A)$) whose fibers are $\V$-connected algebras.
\endcorollary

The Corollary \ref{The fibers of Pierce Sheaf are connected} can be restricted even more, in order to obtain the last result of this paper. The proof is essentially the same of Corollary 14.2 in \cite{CMZ2016}.

\corollary\label{Subdirect product}
Let $\V$ be a variety with RexDFC and CSC. Then, every algebra of $\V$ is a subdirect product of $\V$-connected algebras.
\endcorollary

\subsection*{Acknowledgments}

I would like to thank specially to Prof. Diego Vaggione for his
clarifying comments about the Theory of Central Elements and also for his
useful suggestions about this manuscript.

\refs

\bibitem [Badano2012]{B2012} M. Badano. Variedades con Congruencias Factor Ecuacionalmente Definibles. (PhD Thesis). Universidad Nacional de C\'ordoba, 2012. 
\bibitem [BadanoVaggione2013]{BV2013} M. Badano \& D. Vaggione, Varieties with equationally definable factor congruences, Algebra Universalis. Volume 69, 139 --166, 2013.
\bibitem [BadanoVaggione2017]{BV2017} M. Badano \& D. Vaggione, Varieties with equationally definable factor congruences II, Algebra Universalis. Volume 78, Issue 1 19 -- 42, 2017.
\bibitem [BigelowBurris1990]{BB1990} D.Bigelow \& S. Burris. Boolean algebras of factor congruences, Acta Sci. Math., 54, 11--20, 1990.
\bibitem [Carboni1993]{CW1993}  A. Carboni, S. Lack, R.F.C. Walters. Introduction to extensive and distributive categories. Journal of Pure and Applied Algebra. Volume 84, Issue 2. 145 --158, 1993.
\bibitem[CastiglioniMenniZuluaga2016]{CMZ2016} J.L. Castiglioni, M. Menni \& W.J. Zuluaga Botero. A representation theorem for integral rigs and its applications to
residuated lattices. Journal of Pure and Applied Algebra, 220, 3533--3566, 2016.
\bibitem [FraserHorn1970]{FH1970} G. A. Fraser \& A. Horn, Congruence relations in direct products. Proc. Amer. Math. 26, 390-394, 1970.
\bibitem [Johnstone1982]{J1982} P.T. Johnstone, Stone Spaces, Cambridge Studies in Advanced Mathematics, vol. 3, Cambridge University Press, Cambridge, 1982.
\bibitem [Johnstone2002]{J2002} P.T. Johnstone. Sketches of an elephant: A topos theory compendium. Oxford University Press, 2002.
\bibitem [Lawvere2008]{L2008} F. W. Lawvere, Core varieties, extensivity, and rig geometry. Theory Appl. Categ. 20(14):497-503, 2008.
\bibitem [MacLane1971]{M1971} S. Mac Lane, Categories for the Working Mathematician, Graduate Texts in Mathematics, Springer Verlag, 1971.
\bibitem [MacLaneMoerdijk2012]{MM2012} S. Mac Lane \& I. Moerdijk. Sheaves in geometry and logic: A first introduction to Topos Theory. Springer Science \& Business Media, 2012.
\bibitem [McKenzieMcNultyTaylor1987]{MMT1987} R. McKenzie, G. McNulty \& W. Taylor, Algebras, Lattices, Varieties, Vol. 1, Wadsworth \& BrooksrCole Math. Series, Monterey, CA, 1987.
\bibitem [SanchezVaggione2009]{SV2009} P. Sanchez Terraf \& D. Vaggione, Varieties with definable factor congruences. Trans. Amer. Math. Soc. 361, 50615088, 2009.
\bibitem [Sanchez2010]{S2010} P. Sanchez Terraf, Existentially definable factor congruences, Acta Sci. Math. 76, 49-54, 2010.
\bibitem [Simmons1980]{S1980} H. Simmons, Reticulated rings, J. Algebra 66, 169 -- 192,  1980.
\bibitem [Vaggione1996]{V1996} D. Vaggione, Varieties in which the Pierce Stalks are
Directly Indecomposable. Journal of Algebra. 184, 424--434, 1996.
\bibitem [Vaggione1999] {V1999} D. Vaggione, Central Elements in Varieties with the Fraser -- Horn Property. Advances in Mathematics. Vol. 148, Issue 2, 193 -- 202, 1999. 
\bibitem [Zuluaga2016]{Z2016} W. J. Zuluaga, Representaci\'on por haces de riRigs, PhD Thesis, Universidad Nacional de La Plata, \emph{http://sedici.unlp.edu.ar/handle/10915/54115}, 2016.

\endrefs

\end{document}